\documentclass[12pt, reqno, a4paper]{amsart}
\usepackage{amsmath,mathrsfs}
\usepackage{amssymb}
\usepackage{color}
\usepackage{calligra}
\usepackage{dsfont}
\usepackage{pstricks,pst-node}
\usepackage[latin1]{inputenc}
\usepackage{todonotes}
\usepackage{mathtools} 

\usepackage{pdfsync}

\newcommand\NoBlackBoxes{\global\overfullrule0pt}
\NoBlackBoxes

\usepackage{enumitem}
\setitemize{noitemsep,topsep=0pt,parsep=0pt,partopsep=0pt}

\newcommand{\eps}{\varepsilon}

\newcommand{\N}{\mathbb{N}}

\newcommand{\Z}{\mathbb{Z}}

\renewcommand{\P}{\mathbb{P}}

\newcommand{\Cov}{\mathop{\mathrm{Cov}}\nolimits}

\newcommand{\eee}{{\rm e}}
\newcommand{\dd}{{\rm d}}

\makeatletter
\let\serieslogo@\relax
\let\@setcopyright\relax

\theoremstyle{plain}
\newtheorem{theorem}{Theorem}[section]
\newtheorem{lemma}[theorem]{Lemma}

\newtheorem{proposition}[theorem]{Proposition}

\theoremstyle{definition}

\theoremstyle{remark}

\newtheorem{rem}[theorem]{Remark}

\renewcommand{\P}{{\mathbb{P}}}
\newcommand{\E}{{\mathbb{E}}}
\newcommand{\R}{{\mathbb{R}}}

\newcommand{\C}{\mathbb{C}}

\newcommand{\V}{\mathbb{V}}
\newcommand{\Co}{\mathrm{Cov}}

\renewcommand{\epsilon}{\varepsilon}
\renewcommand{\phi}{\varphi}

\numberwithin{equation}{section}

\begin{document}

\setcounter{page}{1}

\title[Fluctuations for Ising models on Erd\"os-R\'enyi graphs]{Fluctuations for the partition function of Ising models on Erd\"os-R\'enyi random graphs}

\author[Zakhar Kabluchko]{Zakhar Kabluchko}
\address[Zakhar Kabluchko]{Fachbereich Mathematik und Informatik,
Universit\"at M\"unster,
Einsteinstra\ss e 62,
48149 M\"unster,
Germany}

\email[Zakhar Kabluchko]{zakhar.kabluchko@uni-muenster.de}

\author[Matthias L\"owe]{Matthias L\"owe}
\address[Matthias L\"owe]{Fachbereich Mathematik und Informatik,
Universit\"at M\"unster,
Einsteinstra\ss e 62,
48149 M\"unster,
Germany}

\email[Matthias L\"owe]{maloewe@math.uni-muenster.de}

\author[Kristina Schubert]{Kristina Schubert}
\address[Kristina Schubert]{ Fakult\"at f\"ur Mathematik, TU Dortmund, Vogelpothsweg 87, 44227 Dortmund,
Germany}

\email[Kristina Schubert]{kristina.schubert@tu-dortmund.de}


\date{\today}



\newcommand{\wlim}{\mathop{\hbox{\rm w-lim}}}
\newcommand{\na}{{\mathbb N}}
\newcommand{\re}{{\mathbb R}}

\newcommand{\vep}{\varepsilon}

\begin{abstract}
We analyze Ising/Curie-Weiss models on the Erd\H{o}s-R\'enyi graph with $N$ vertices and edge probability $p=p(N)$ that were introduced by Bovier and Gayrard [J.\ Statist.\ Phys., 72(3-4):643--664, 1993] and investigated in \cite{KLS19a} and \cite{KLS19b}. We prove  Central Limit Theorems for the partition function of the model and -- at other decay regimes of $p(N)$  -- for the logarithmic partition function. We find critical regimes for $p(N)$ at which the behavior of the fluctuations of the partition function changes.
\end{abstract}

\maketitle

\section{Introduction}
In this note we add another step to our analysis of Ising models on the Erd\H{o}s-R\'enyi random graph. These models are defined on a realization of a finite directed  random graph $G=(V,E)=(\{1,\ldots, N\},E)$ by the Hamiltonian or energy function
\begin{equation}\label{hamil}
H(\sigma)= H_N(\sigma)= - \frac 1 {2Np} \sum_{i,j=1}^N \sigma_i \sigma_j \vep_{i,j}, \qquad \sigma \in \{-1,+1\}^N.
\end{equation}
Here $N$ is the number of vertices of the graph, $p=p(N)$ is the probability that a directed edge $(i,j)$ is present in $E$ and $\vep_{i,j}$ is the indicator for this event. Moreover, we assume these indicator variables to be independent.
With $H$ we associate a Gibbs measure $\mu_\beta$ at inverse temperature $\beta>0$ defined as
\begin{equation}\label{gibbs}
\mu_\beta (\sigma)\coloneqq\frac 1 {Z_{N,\beta}} \exp(-\beta H(\sigma)), \qquad \sigma \in \{-1,+1\}^N.
\end{equation}
The quantity
\begin{equation}\label{partition}
Z_{N,\beta} \coloneqq Z_N(\beta)\coloneqq\sum_{\sigma \in \{-1,+1\}^N}\exp(-\beta H(\sigma))
\end{equation}
is called the partition function of the model. It encodes much of the interesting information about the system and will feature as the main character in the present note.
The limit
\begin{equation}\label{free_energy}
\lim_{N \to \infty} -\frac 1 {\beta N} \log Z_{N,\beta},
\end{equation}
if it exists, is called the free energy per site or particle.

The easiest nontrivial case of the above situation is the Curie-Weiss model, i.e.\ the case where $p=1$. Here the Hamiltonian can be rewritten as
$$
H(\sigma) = - \frac N 2 \left(\frac{\sum_{i=1}^N \sigma_i}{N}\right)^2 =: - \frac N 2 m_N(\sigma)^2
$$
i.e.\ it is a function of the magnetization per particle $m_N(\sigma)= \frac 1N \sum_{i=1}^N \sigma_i$,
which makes it accessible, among others, to large deviation theory.

The Curie-Weiss model was studied in great detail. Convergence results can be found e.g.~in \cite{Ellis_Newman_78a}, \cite{Ellis_Newman_78b}, \cite{EiseleEllis-MultiplePhaseTransitionsInGCW}, or the monograph \cite{Ellis-EntropyLargeDeviationsAndStatisticalMechanics}. One major result is a phase transition at $\beta =1$. While in the high temperature regime $\beta \le 1$ the magnetization $m_N$ converges to 0 as the system size $N$ goes to infinity, its distribution under the Gibbs measure converges to the mixture $\frac 12 (\delta_{m^+(\beta)}+\delta_{m^-(\beta)})$, if $\beta >1$ (the low temperature regime). Here $\delta_x$ denotes the Dirac measure in a point $x$, $m^+(\beta)$ is the largest solution of
$$
z = \tanh(\beta z),
$$
and $m^-(\beta)=-m^+(\beta)$.

Moreover, in \cite{Ellis_Newman_78b}, \cite{Ellis-EntropyLargeDeviationsAndStatisticalMechanics}, \cite{EL10}, \cite{Chatterjee_Shao} it was shown that this phase transition in the Curie-Weiss model is also visible on the level of fluctuations. While for $\beta <1$ the rescaled magnetization
$\sqrt N m_N$ converges in distribution to a centered normal random variable with variance $\frac 1{1-\beta}$, for $\beta =1$ one has to scale differently. Here one obtains that $\sqrt[4] N m_N$ converges in distribution to a non-normal random variable
with Lebesgue density proportional to $\exp(-\frac 1 {12} x^4)$.

In \cite{BG93b} it was shown that the law of large number type results for $m_N$ still hold true for the above
Ising model on a random graph, as long as $Np \to \infty$ (which implies that almost surely the graph has a giant component that contains almost all vertices). Indeed, in this case almost surely with respect to the probability measure that describes the random graph, the quantity $m_N$ behaves as in the Curie-Weiss model. This situation has to be contrasted to the results by Dembo and Montanari in \cite{Dembo_Montanari_2010a} and \cite{Dembo_Montanari_2010b} as well as Giardina and van der Hofstad with coauthors in \cite{van_der_Hofstad_et_al_2010},
\cite{van_der_Hofstad_et_al_2014}, \cite{van_der_Hofstad_et_al_2015}, \cite{van_der_Hofstad_et_al_2015b}, and \cite{van_der_Hofstad_et_al_2015c}, who treat the difficult case of locally tree-like, i.e.\ sparse, random graph models, and analyze thermodynamic quantities there.

In our papers \cite{KLS19a} and \cite{KLS19b} we analyzed the fluctuations of $m_N$ in the Ising model on an Erd\H{o}s-R\'enyi graph given by \eqref{gibbs} in the regime $Np\to\infty$. We were able to show that for $\beta <1$ the following holds true:
Consider the following {\it random} probability measures on $\R$
$$
L_N \coloneqq \frac 1{Z_N(\beta)} \sum_{\sigma\in \{-1,+1\}^N} \eee^{-\beta H(\sigma)}  \delta_{\frac 1 {\sqrt N} \sum_{i=1}^N \sigma_i}.
$$
Then $L_N$, viewed as a random element in the space of probability measures on $\R$ endowed with the weak topology,  converges in probability to the normal distribution with mean 0 and variance $\frac{1}{1-\beta}$.
For $\beta=1$ we obtained that
$$
\overline L_N \coloneqq \frac 1{Z_N(\beta)} \sum_{\sigma\in \{-1,+1\}^N} \eee^{-\beta H(\sigma)}  \delta_{\frac 1 {N^{3/4}} \sum_{i=1}^N \sigma_i}
$$
converges in probability to a non-normal random variable with Lebesgue density proportional to $\exp(-\frac 1 {12} x^4)$, if $p^4 N^3 \to \infty$, while for regimes with smaller $p$ an appropriate modification of $\overline L_N$ converges to a normal distribution.

A key tool in \cite{KLS19a} was that for $p$ so large that $p^3N^2 \to \infty$ and any bounded continuous function $g \in \mathcal{C}^b(\R), g\ge 0$ the generalized partition function
\begin{equation}\label{ZN(g)}
Z_N(\beta, g)\coloneqq\sum_{\sigma \in \{-1, +1\}^N} \eee^{-\beta H(\sigma)} g\left( \frac{\sum_{i=1}^N \sigma_i}{\sqrt N} \right)
\end{equation}
satisfies
\begin{equation}\label{eq:convergence}
\frac{Z_N(\beta, g)}{\E Z_N(\beta,g)} \to 1
\end{equation}
in $L^2$ and, hence, in probability.
Here, $\E$ denotes expectation with respect to the probability measure $\P$, i.e.~the randomness generated by $(\eps_{i,j})_{i,j=1}^N$.

Applying this to $g \equiv 1$ proves a Law of Large Numbers for the partition function:
\begin{equation}\label{eq:convergence2}
\frac{Z_N(\beta)}{\E Z_N(\beta)} \to 1
\end{equation}
in $L^2$ and, hence, in probability. This holds in the regime $p^3 N^2 \to \infty$. Moreover, in the same regime, the denominator of the fraction in \eqref{eq:convergence2} was shown to satisfy
$$
\E Z_N (\beta)\sim \frac{2^N}{\sqrt{1-\beta}} e^{\frac{(1-p)\beta^2}{8p}}.
$$
Here, 
for two sequences $(a_N), (b_N)$ we write $a_N \sim b_N$, if and only if their quotient converges to $1$.

Seeing a Law of Large Numbers like \eqref{eq:convergence2} one is automatically tempted to ask for a Central Limit Theorem (CLT, for short) for the fluctuations around it. The aim of the present note is to give an answer to this question. Similar fluctuation results for the partition function of disordered systems like SK-models or the REM were shown in \cite{ALR}, \cite{CN}, and \cite{BovierKurkovaLoewe-FluctuationsOfFreeEnergyREM}.

More precisely, for the rest of the paper we will assume that $\beta <1$ and we will show Gaussian fluctuations for the partition function $Z_N(\beta)$ as long as $p^3 N^2 \to \infty$ and Gaussian fluctuations for the log-partition function $\log Z_N(\beta)$, if $p^3 N^2 \to c \in [0, \infty)$. This also shows that the model indeed undergoes a phase transition in the parameter regime when $p^3 N^2 \to c \in (0, \infty)$. In \cite{KLS19a} and \cite{KLS19b} this regime was the regime where we had to change the techniques to prove a CLT for the magnetization, however, the results did not indicate a phase transition.
We will show the following theorems.

\begin{theorem}\label{CLT_ZN_1}
Consider the Ising model on an Erd\H{o}s-R\'enyi random graph. Define $Z_N(\beta)$ as in \eqref{partition}.  Suppose that $0< \beta <1$, $p^2N  \to \infty$, and $(1-p) Np \to \infty$ (the latter condition is fulfilled if $p(N)$ is bounded away from $1$). Then,
\begin{equation}
\sqrt {\frac{pN}{1-p}} \left( \frac{Z_N(\beta)}{\E Z_N(\beta)} -1\right) \to \mathcal{N}\left(0,\frac {\beta^2}{4} \right)
\end{equation}
in distribution.
\end{theorem}

The previous theorem and its proof immediately raise the question what happens, if the probability $p$ of retaining an edge in the Erd\H{o}s-R\'enyi graph gets even smaller, such that $p^2 N \to \infty$ no longer holds true. The following theorem gives an answer to this problem.

\begin{theorem}\label{CLT_ZN_2}
Consider the Ising model on an Erd\H{o}s-R\'enyi random graph. Define $Z_N(\beta)$ as in \eqref{partition}. Assume that $0< \beta <1$.
\begin{itemize}
\item
If $p \sqrt N  \to c \in (0,\infty)$, then
\begin{equation}
N^{1/4} \left( \frac{Z_N(\beta)}{\E Z_N(\beta)} -1\right) \to \mathcal{N}\left(0,\frac{\beta^2}{4c^2}+\frac{\beta^4}{64c^3}\right)
\end{equation}
in distribution.
\item
If $p^2N  \to 0$ but $p^3 N^2 \to \infty$, then
\begin{equation}
N p^{3/2} \left( \frac{Z_N(\beta)}{\E Z_N(\beta)} -1\right) \to \mathcal{N}\left(0,\frac{\beta^4}{64}\right)
\end{equation}
in distribution.
\end{itemize}
\end{theorem}

If $p^3 N^2 \to c \in (0,\infty)$, we obtain normal fluctuations for the log-partition function $\log Z_N(\beta)$  and log-normal fluctuations for the partition function $Z_N(\beta)$.

\begin{theorem}\label{CLT_ZN_3}
If $p^3 N^2 \to c \in (0,\infty)$ and $0<\beta <1$, then
\begin{equation*}
\log Z_N(\beta)- N\log 2 - \frac {\beta^2(1-p)}{8p} + \frac 12 \log (1-\beta) + \frac{\beta^4}{192 c} \to \mathcal{N}\left(0,\frac{\beta^4}{64c}\right)
\end{equation*}
in distribution. In other words,
$$
\sqrt{1-\beta} \, \frac{Z_N(\beta)}{2^N e^{\frac {\beta^2(1-p)}{8p} -\frac{\beta^4}{192 c}}} \to e^{\frac{\beta^2\xi}{8\sqrt c}}
$$
in distribution, where $\xi$ is a standard normal random variable.
\end{theorem}

Finally, if $p^3 N^2 \to 0$, then the log-partition function $\log Z_N(\beta)$ still has normal fluctuations, but the size of these fluctuations, namely $1/(Np^{3/2})$, diverges as $n\to\infty$ and hence there is no affine normalization of the form  $a_N(Z_N(\beta)-b_N)$ for the partition function itself which converges to a non-degenerate limit distribution.

\begin{theorem}\label{CLT_ZN_4}
If $p^3 N^2 \to 0$ and $0<\beta <1$, then
\begin{equation*}
Np^{3/2}\left(\log Z_N(\beta)- N\log 2 - N^2p \log \cosh\left(\frac{\beta}{2Np}\right)  \right)
\to \mathcal{N}\left(0,\frac{\beta^4}{64}\right)
\end{equation*}
in distribution.
\end{theorem}

The article is organized in the following way. In the next section we will state some asymptotic expansion results needed in the main proofs. These are given in Section~3 (for Theorem~\ref{CLT_ZN_1}), Section~4 (for Theorem~\ref{CLT_ZN_2}), and Section~5 (for Theorems~\ref{CLT_ZN_3} and~\ref{CLT_ZN_4}).

\section{Technical Preparations}
In this section we will recall some results needed in the next sections.
For the rest of this note, given spin configurations $\sigma, \tau \in\{-1,+1\}^N$ let us define
\begin{equation}
|\sigma| \coloneqq \sum_{i=1}^N \sigma_i
\end{equation}
as well as
\begin{equation}
|\sigma\tau| \coloneqq \sum_{i=1}^N \sigma_i \tau_i.
\end{equation}

Next, for arbitrary complex variables $z$ and $p$ and an integer $m$ let
us define the function
\begin{equation}\label{eq:centralfunc}
F_m(x,p, z) \coloneqq \log \left(1 - p + p\eee^{xz-m(1-2p)\frac{z^2}{2}}\right).
\end{equation}

Let us stress that in this part of the paper we do not require that $p$ is a probability.
Next we will state some power series expansion involving $F_m(x,p, z)$. Here, the variables $x$ and $m$ will be fixed and the function will be expanded in the $p$ and $z$ variables around the origin $(0,0)$.
For any fixed values of $x$ and $m$, and for $|p|< 2$ and $|z|< z_0$ with sufficiently small $z_0>0$, we have
$$
|p\eee^{xz-m(1-2p)\frac{z^2}{2}}-p| < 1.
$$
Thus, for fixed values of $x$ and $m$, the function $F_m(x,p,z)$ is an analytic function of two complex variables $p$ and $z$ on the polydisc domain
$$
\mathcal D = \{(p,z)\in\C^2\colon |p| < 2, |z|<z_0\}.
$$
Therefore, it has a power series expansion which converges uniformly and absolutely on compact subsets of this domain. In particular, we may re-arrange and re-group the terms arbitrarily without changing the sum.
In the sequel we will use the first terms of the power series expansion of some expressions involving $F_m(x,p,z)$. These are
given in the following

\begin{lemma}\label{taylor}
For linear combinations involving $F_1$ we have  power series expansions around the point $(0,0)$ of the form
\begin{align}
\frac {F_1(1, p, z) + F_1(-1, p, z)}2 = &\frac 12 p^2z^2 - \frac 1{12} pz^4 +  \mathcal{O}(p^2z^4) +  \mathcal{O}(pz^6),\label{f1a}\\
\frac {F_1(1, p, z) - F_1(-1, p, z)}2 = &pz +  \mathcal{O}(pz^3).\label{f1b}
\end{align}
For linear combinations involving $F_2$ we have power series expansions around the point $(0,0)$ of the form
\begin{align}
\frac {F_1(2, p, z) + F_1(-2, p, z)+2F_1(-2, p, z)}4
&=
p^2z^2 - \frac 1{6} pz^4 +  \mathcal{O}(p^2z^4) +  \mathcal{O}(pz^6),\label{f2a}\\
\frac {F_1(2, p, z) + F_1(-2, p, z)-2F_1(-2, p, z)}4
&=
p(1-p)z^2 - \frac 2{3} pz^4 +  \mathcal{O}(p^2z^4) +  \mathcal{O}(pz^6),\label{f2b}\\
\frac {F_1(2, p, z) + F_1(-2, p, z)}4
&=
pz + \mathcal O(pz^3). \label{f2c}
\end{align}
\end{lemma}
\begin{proof}
An easy Taylor expansion in the $p$ and $z$ variables gives the above formulas.
\end{proof}

The following Lemmas~\ref{lem:expect} and~\ref{lem:lemma_cov} were shown in \cite{KLS19a}, Theorem 3, equations (4.7)--(4.10), and Lemmas~2 and~3.
\begin{lemma}\label{lem:expect}
Assume that $p^3N^2 \to \infty$ and set
\begin{equation}\label{eq:A_N_beta_def}
A_{N}(\beta) \coloneqq -\frac{\beta^2}{8}+ N^2p \left(\cosh\left(\frac{\beta}{2Np}\right)-1\right) = -\frac{\beta^2}{8}+\frac{\beta^2}{8p}+o(1).
\end{equation}
Then we have for all $\sigma \in\{\pm 1\}^N$
\begin{align}
\E [\eee^{-\beta H(\sigma)}]
&= 
\eee^{A_N(\beta) + \frac{\beta}{2} \frac {|\sigma|^2} N + \frac 1 {N^2p^2}   \left(C_{N,1} + C_{N,2}\frac {|\sigma|^2}N\right)}. \label{eq:exp_computation1}
\end{align}
Here the sequences $(C_{N,1})_{N\in\N}$ and $(C_{N,2})_{N\in\N}$ do not depend on $\sigma\in \{-1, +1\}^N$ and are bounded.
For the expected partition function this yields
\begin{equation}\label{eq EZNg}
\E Z_N (\beta)\sim \eee^{\frac{(1-p)\beta^2}{8p}} \frac {2^N} {\sqrt{1-\beta}}.
\end{equation}
\end{lemma}

\begin{lemma}\label{lem:lemma_cov}
Assume that $p^3N^2\to \infty$ and define $A_N(\beta)$ as in Lemma \ref{lem:expect}, which yields
\begin{equation}\label{eq:A_N_beta_asympt}
2A_{N}(\beta) = -\frac{\beta^2}{4}+ 2N^2p \left(\cosh\left(\frac{\beta}{2Np}\right)-1\right) = -\frac{\beta^2}{4}+\frac{\beta^2}{4p}+
\mathcal{O}\left(\frac 1 {N^2 p^3}\right).
\end{equation}
Moreover, let
\begin{equation}\label{eq:B_N_beta_asympt}
B_{N}(\beta) \coloneqq -\frac{\beta^2}{4}+\frac{N^2p}{2} \left(\cosh\left(\frac{\beta}{Np}\right)-1\right)
=
-\frac{\beta^2}{4}+\frac{\beta^2}{4p}+\mathcal{O}\left(\frac 1 {N^2 p^3}\right).
\end{equation}
Then
\begin{align}
\E [\eee^{-\beta H(\sigma)}] \E [\eee^{-\beta H(\tau)}]
&= 
\eee^{2A_N(\beta) + \frac{\beta}{2} \frac {|\sigma|^2+|\tau|^2} N + \frac 1 {N^2p^2}   \left(2C_{N,1} + C_{N,2}\frac {|\sigma|^2 + |\tau|^2}N\right)}, \label{eq:cov_computation1}\\
\E [\eee^{-\beta H(\sigma)} \eee^{-\beta H(\tau)}]
&=
\eee^{
\left(B_N(\beta) + \frac{C_{N,3}}{N^2p^2}\right) \left(1+\frac{|\sigma\tau|^2}{N^2}\right)
+
\frac{|\sigma|^2+|\tau|^2}N \left(\frac \beta 2 + \frac{C_{N,4}}{N^2p^2} \right)
}.\label{eq:cov_computation2}
\end{align}
Here the sequences $(C_{N,1})_{N\in\N}$ and $(C_{N,2})_{N\in\N}$ are the same  as in Lemma \ref{lem:expect} and $(C_{N,3})_{N\in\N}$, and $(C_{N,4})_{N\in\N}$ also do not depend on $\sigma, \tau \in \{-1, +1\}^N$ and are bounded.

Hence
\begin{multline*}
\Co(\eee^{-\beta H(\sigma)}, \eee^{-\beta H(\tau)})
=
\E [\eee^{-\beta H(\sigma)}] \E [\eee^{-\beta H(\tau)}]
 \left(\eee^{\frac{|\sigma\tau|^2}{N^2}\frac{\beta^2(1-p)}{4p}+\mathcal{O}\left(\frac 1 {N^2 p^3} + \frac {|\sigma|^2+|\tau|^2} {N^3p^2}\right)}-1\right).
\end{multline*}
\end{lemma}
We will also need the following expansion which is easily checked.
\begin{lemma}\label{lem:expan_exp}
For real $x$ in a neighbourhood of $0$ and $|p| \le 2$ we have
$$
\frac{\eee^x}{1-p+p \eee^{x}}-1=(1-p)x +\mathcal{O}(x^2).
$$
\end{lemma}

\section{Proof of Theorem \ref{CLT_ZN_1}}\label{sec:proof_CLT_ZN_1}

\subsection{Method of proof}
The proof of Theorem \ref{CLT_ZN_1} relies on an approximation of $Z_N(\beta)$ by terms that are more easily seen to fulfill a CLT.
Let us motivate this.
The idea of the proof is to take a Hermite expansion of the exponential in the partition function and to show that the first few terms are dominant. For these,  in turn, it will be sufficiently simple to derive a CLT. Thus, let $\mathit{He}_n(x)$ be the $n$'th probabilistic Hermite polynomial, i.e.
$$
\mathit{He}_n(x)\coloneqq 
(-1)^n \eee^{x^2/2}\frac {d^n}{dx^n} \eee^{-x^2/2}.
$$
Recall that for a standard normal random variable $\xi$ the Hermite expansion of a function $f(\xi)$ is given by
$$
f(\xi)= \sum_{n=0}^\infty \mathit{He}_n(\xi) \hat f_n,
$$
where
$
\hat f_n \coloneqq \frac 1 {n!}  \E (f(\xi) \mathit{He}_n(\xi)).
$
To ``see'' the Gaussian random variable in our situation, notice that for all $\sigma\in \{-1,+1\}^N$,
$$
\E H_N(\sigma) = -\frac 1 {2N} |\sigma|^2 \qquad \mbox{and } \qquad  \V H_N(\sigma)= \frac{1-p}{4p},
$$
such that we may be tempted to approximate $H_N(\sigma)$ in the following way:
$$
H_N(\sigma) \approx -\frac 1 {2N} |\sigma|^2 - \sqrt{\frac{1-p}{4p}} \xi, \qquad \xi \sim \mathcal{N}(0,1),
$$
and
$$
\eee^{-\beta H_N(\sigma) }\approx \eee^{\beta (\frac 1 {2N} |\sigma|^2 + \sqrt{\frac{1-p}{4p}} \xi)} =: \eee^{a+b \xi}.
$$
Since $\mathit{He}_0(x)=1$ and $\mathit{He}_1(x)=x$ we obtain as the first coefficients of the Hermite expansion of $f(\xi)\coloneqq\eee^{a+b \xi}$:
$$
\hat f _0 = \E \eee^{a+b \xi}= \eee^{a + \frac{b^2}2} \qquad \mbox{and } \qquad \hat f_1 =\E\left( \eee^{a+b \xi}\xi \right)= b  \eee^{a + \frac{b^2}2}
$$
and thus (hopefully)
\begin{eqnarray*}
\eee^{-\beta H_N(\sigma)}&\approx& \eee^{a + \frac{b^2}2}(1+ b\xi) =  \eee^{a + \frac{b^2}2}(a+ b\xi+1-a)\\
&\approx &   \eee^{a + \frac{b^2}2}(-\beta H_N(\sigma) +1 - \frac{\beta}{2N} |\sigma|^2) \\
&=& \eee^{\frac{\beta}{2N}|\sigma|^2 +\frac{\beta^2(1-p)}{8p}}\left[1-\beta H_N(\sigma)- \frac{\beta}{2N} |\sigma|^2\right].
\end{eqnarray*}
As the first factor on the right is approximately $\E \eee^{-\beta H_N(\sigma)}$ (cf. e.g.\, Lemma 2 in \cite{KLS19a}) we will work with the approximation
\begin{equation}\label{eq:approx}
\eee^{-\beta H_N(\sigma)}\approx \E (\eee^{-\beta H_N(\sigma)})\left[1-\beta H_N(\sigma)- \frac{\beta}{2N} |\sigma|^2\right].
\end{equation}
The hope is that the two expressions are close together in $L^2$ and that for the term on the right hand side we can obtain a Central Limit Theorem.

\subsection{Expectations and covariances}
To justify these ideas put
\begin{align}\label{Xsigma}
X(\sigma)&\coloneqq  \eee^{-\beta H(\sigma)}- \E (\eee^{-\beta H(\sigma)})\left[1-\beta H(\sigma)- \frac{\beta}{2N} |\sigma|^2\right]\nonumber \\
&= \eee^{-\beta H(\sigma)}- \E \eee^{-\beta H(\sigma)}+\E (\eee^{-\beta H(\sigma)})\left[\beta H(\sigma)+ \frac{\beta}{2N} |\sigma|^2\right],
\end{align}
where we write $H(\sigma):= H_N(\sigma)$ to simplify the notation. We begin with
\begin{lemma}
For all $\sigma \in\{-1,+1\}^N$ we have $\E X(\sigma)=0$.
\end{lemma}
\begin{proof}
This is easily checked using the second line in \eqref{Xsigma}.
\end{proof}

Our aim is to check that $\sum_\sigma X(\sigma)$ has ``small'' variance, which would justify the approximation in~\eqref{eq:approx}.
Now,
$$
\V\sum_\sigma X(\sigma)= \sum_{\sigma, \tau} \Co(X(\sigma),X(\tau))
$$
and (bearing in mind that
$\eee^{-\beta H(\sigma)}- \E \eee^{-\beta H(\sigma)}$ and $\beta H(\sigma)+ \frac{\beta}{2N} |\sigma|^2 $ are centered),
\begin{eqnarray}\label{Cov1a}
\Co(X(\sigma),X(\tau))&= &\Co(\exp(-\beta H(\sigma), \exp(-\beta H(\tau)) \nonumber \\
&&+\E \eee^{-\beta H(\sigma)}\E \eee^{-\beta H(\tau)} \Co(\beta H(\sigma), \beta H(\tau))\nonumber \\
&&+ \E \eee^{-\beta H(\sigma)} \Co(\eee^{-\beta H(\tau)}, \beta H(\sigma))\nonumber \\
&&+\E \eee^{-\beta H(\tau)} \Co(\eee^{-\beta H(\sigma)}, \beta H(\tau)).
\end{eqnarray}
Let us compute each of the four terms in \eqref{Cov1a} separately.
\begin{lemma}\label{lem:covA}
For all $p$ satisfying $N^2 p^3 \to \infty$ and all $\sigma, \tau \in \{-1,+1\}^N$ we have
$$
\Co(\eee^{-\beta H(\sigma)}, \eee^{-\beta H(\tau)})
=
\E [\eee^{-\beta H(\sigma)}] \E [\eee^{-\beta H(\tau)}]
 \left(\eee^{\frac{|\sigma\tau|^2}{N^2}\frac{\beta^2(1-p)}{4p}+\mathcal{O}\left(\frac 1 {N^2 p^3} + \frac {|\sigma|^2+|\tau|^2} {N^3p^2}\right)}-1\right).
$$
\end{lemma}
\begin{proof}
This is just a restatement of Lemma \ref{lem:lemma_cov}.
\end{proof}

For the second term in \eqref{Cov1a} we obtain
\begin{lemma}\label{lem:covB}
For all $\sigma, \tau \in \{-1,+1\}^N$ and all $p$ we have
$$
\Co(\beta H(\sigma), \beta H(\tau))=\frac{\beta^2(1-p)}{4N^2p}|\sigma \tau|^2.
$$
\end{lemma}
\begin{proof}
This follows from
\begin{eqnarray*}
\Co(\beta H(\sigma), \beta H(\tau))&=& \beta^2 \Co\left(\frac 1 {2Np} \sum_{i,j} \vep_{i,j} \sigma_i \sigma_j, \frac 1 {2Np} \sum_{i,j} \vep_{i,j} \tau_i \tau_j\right)\\
&=&\beta^2\left( \frac 1{4N^2p^2}\E(\sum_{(i,j), (i',j')} \vep_{i,j}  \vep_{i',j'}  \sigma_i \sigma_j \tau_{i'} \tau_{j'})-\frac 1{4 N^2}|\sigma|^2 |\tau|^2\right)\\
&=&\beta^2\left(\frac 1{4N^2} \sum_{(i,j) \neq (i',j')}  \sigma_i \sigma_j \tau_{i'} \tau_{j'}+\frac 1{4N^2p}\sum_{(i,j)}  \sigma_i \sigma_j \tau_i \tau_j -\frac 1{4 N^2}|\sigma|^2 |\tau|^2\right)\\
&=&\frac{\beta^2(1-p)}{4N^2p}|\sigma \tau|^2.
\end{eqnarray*}
\end{proof}
Hence we obtain
$$
\E \eee^{-\beta H(\sigma)}\E \eee^{-\beta H(\tau)} \Co(\beta H(\sigma), \beta H(\tau))= \E \eee^{-\beta H(\sigma)}\E \eee^{-\beta H(\tau)}\frac{\beta^2(1-p)}{4N^2p}|\sigma \tau|^2.
$$
For the third and the fourth term in \eqref{Cov1a} we need to compute a covariance of the form $ \Co(\eee^{-\beta H(\sigma)}, \beta H(\tau))$:
\begin{lemma}\label{lem:covC}
For all $\sigma, \tau \in \{-1,+1\}^N$ and all $p$ we have
$$
\Co(\eee^{-\beta H(\sigma)}, \beta H(\tau))=
-\E \eee^{-\beta H(\sigma)}\left(\frac{\beta^2 (1-p)}{4N^2 p}  | \sigma \tau|^2+ \mathcal{O}\left(\frac {|\tau|^2}{N^2 p^2}\right)\right).
$$
\end{lemma}
\begin{proof}
We begin with
\begin{eqnarray*}
\E[\eee^{-\beta H(\sigma)} \beta H(\tau)]&=& -\E\left[ \prod_{i,j} \exp \left(\frac{\beta}{2Np} \sigma_i \sigma_j \vep_{i,j}\right) \sum_{k,l}  \frac{\beta}{2Np} \tau_k \tau_l \vep_{k,l}\right]\\
&=& -\sum_{k,l} \E \left[ \frac{\beta}{2Np} \tau_k \tau_l \vep_{k,l}\prod_{i,j} \exp \left(\frac{\beta}{2Np} \sigma_i \sigma_j \vep_{i,j}\right) \right]\\
&=&  -\sum_{k,l} \E \left[ \frac{\beta}{2Np} \tau_k \tau_l \vep_{k,l}\eee^{\frac{\beta}{2Np} \sigma_k \sigma_l \vep_{k,l}}\prod_{(i,j) \neq (k,l)} \eee^{\frac{\beta}{2Np} \sigma_i \sigma_j \vep_{i,j}} \right]
\end{eqnarray*}
by differentiating the cases $(i,j)=(k,l)$ and $(i,j) \neq (k,l)$. Exploiting the independence of the $\vep_{i,j}$'s we obtain:
\begin{eqnarray*}
\E[\eee^{-\beta H(\sigma)} \beta H(\tau)]
&=&
 -\sum_{k,l} \frac{\beta}{2N}\tau_k \tau_l \eee^{\frac{\beta}{2Np} \sigma_k \sigma_l}
\left(\prod_{i,j}(1-p+p \eee^{\frac{\beta}{2Np} \sigma_i \sigma_j})\right) \frac 1 {1-p+p\eee^{\frac{\beta}{2Np} \sigma_k \sigma_l}}\\
&=& -\E \eee^{-\beta H(\sigma)} \times  \sum_{k,l} \frac{\frac{\beta}{2N}\tau_k \tau_l \eee^{\frac{\beta}{2Np} \sigma_k \sigma_l}}{1-p+pe^{\frac{\beta}{2Np} \sigma_k \sigma_l}}
\end{eqnarray*}
as the product in the first line equals $\E \eee^{-\beta H(\sigma)}$.
We thus arrive at
\begin{eqnarray*}
\Co(\eee^{-\beta H(\sigma)}, \beta H(\tau))&=&
 - \E \eee^{-\beta H(\sigma)} \times \sum_{k,l}\left( \frac{\frac{\beta}{2N}\tau_k \tau_l \exp (\frac{\beta}{2Np} \sigma_k \sigma_l)}{1-p+p\eee^{\frac{\beta}{2Np} \sigma_k \sigma_l}}\right)+\E \eee^{-\beta H(\sigma)}\frac{\beta}{2N}|\tau|^2\\
&=& - \E \eee^{-\beta H(\sigma)} \times \sum_{k,l}\left( \frac{\frac{\beta}{2N}\tau_k \tau_l \exp (\frac{\beta}{2Np} \sigma_k \sigma_l)}{1-p+p\eee^{\frac{\beta}{2Np} \sigma_k \sigma_l}}-\frac{\beta}{2N}\tau_k \tau_l \right)\\
&=&  - \E \eee^{-\beta H(\sigma)}\frac{\beta}{2N}  \sum_{k,l}\tau_k \tau_l\left( \frac{ \eee^{\frac{\beta}{2Np} \sigma_k \sigma_l}}{1-p+p\eee^{\frac{\beta}{2Np} \sigma_k \sigma_l}}-1\right).
\end{eqnarray*}
Applying Lemma \ref{lem:expan_exp} to  $\frac{ \eee^{\frac{\beta}{2Np} \sigma_k \sigma_l}}{1-p+p\eee^{\frac{\beta}{2Np} \sigma_k \sigma_l}}-1$ we see that
\begin{eqnarray*}
\Co(\eee^{-\beta H(\sigma)}, \beta H(\tau))
&=&   -\E \eee^{-\beta H(\sigma)}\frac{\beta}{2N}  \sum_{k,l}\tau_k \tau_l (1-p) \frac{\beta}{2Np}\left(\sigma_k \sigma_l +\mathcal{O}\left(\frac 1 {N p}\right)\right)\\
&=&   -\E \eee^{-\beta H(\sigma)}\frac{\beta^2 (1-p)}{4N^2 p}  \sum_{k,l}\tau_k \tau_l \left(\sigma_k \sigma_l +\mathcal{O}\left(\frac 1 {N p}\right)\right)\\
&=& -\E \eee^{-\beta H(\sigma)}\left(\frac{\beta^2 (1-p)}{4N^2 p}  | \sigma \tau|^2+ \mathcal{O}\left(\frac {|\tau|^2}{N^3 p^2}\right)\right)
\end{eqnarray*}
where the $\mathcal{O}$-terms do not depend on $\sigma$.
\end{proof}
Therefore we obtain for the third and the fourth term in \eqref{Cov1a}
$$
\E \eee^{-\beta H(\sigma)} \Co(\eee^{-\beta H(\tau)}, \beta H(\sigma)) =
-\E \eee^{-\beta H(\sigma)}\E \eee^{-\beta H(\tau)}\left(\frac{\beta^2 (1-p)}{4N^2 p}  | \sigma \tau|^2+ \mathcal{O}\left(\frac {|\sigma|^2}{N^3 p^2}\right)\right)
$$
as well as
$$
\E \eee^{-\beta H(\tau)} \Co(\eee^{-\beta H(\sigma)}, \beta H(\tau)) =
-\E \eee^{-\beta H(\sigma)}\E \eee^{-\beta H(\tau)}\left(\frac{\beta^2 (1-p)}{4N^2 p}  | \sigma \tau|^2+ \mathcal{O}\left(\frac {|\tau|^2}{N^3 p^2}\right)\right).
$$

\subsection{Estimate for the variance}
Putting the results of Lemmas \ref{lem:covA}--\ref{lem:covC} together we obtain the following estimate for the covariances in \eqref{Cov1a}:
\begin{multline*}
\Co(X(\sigma),X(\tau))= \E \eee^{-\beta H(\sigma)}\E \eee^{-\beta H(\tau)}
\left(\eee^{\frac{|\sigma\tau|^2}{N^2}\frac{\beta^2(1-p)}{4p}+
\mathcal{O}\left(\frac 1 {N^2 p^3} + \frac {|\sigma|^2+|\tau|^2} {N^3p^2}\right)}-1\right.
\\
\left.+\frac{\beta^2 (1-p)}{4N^2 p}  | \sigma \tau|^2-2\left(\frac{\beta^2 (1-p)}{4N^2 p}  | \sigma \tau|^2\right)+ \mathcal{O}\left(\frac {|\sigma|^2}{N^3 p^2}\right)+\mathcal{O}\left(\frac {|\tau|^2}{N^3 p^2}\right)\right).
\end{multline*}
Thus we have
\begin{eqnarray*}
&&\V\left(\sum_\sigma X(\sigma)\right)= \sum_{\sigma, \tau}\E \eee^{-\beta H(\sigma)}\E \eee^{-\beta H(\tau)}
\left(\eee^{\frac{|\sigma\tau|^2}{N^2}\frac{\beta^2(1-p)}{4p}+
\mathcal{O}\left(\frac 1 {N^2 p^3} + \frac {|\sigma|^2+|\tau|^2} {N^3p^2}\right)}-1\right.
\\
&& \left.+\frac{\beta^2 (1-p)}{4N^2 p}  | \sigma \tau|^2-2\left(\frac{\beta^2 (1-p)}{4N^2 p}  | \sigma \tau|^2\right)+ \mathcal{O}\left(\frac {|\sigma|^2}{N^3 p^2}\right)+\mathcal{O}\left(\frac {|\tau|^2}{N^3 p^2}\right)\right).
\end{eqnarray*}
We shall use this to prove the following key result justifying the approximation in~\eqref{eq:approx}:
\begin{lemma}\label{lem:varA}
Assume that $N^2 p^3 \to \infty$. Then,
\begin{equation}\label{estimate_var_sum_X}
\V\left(\sum_\sigma X(\sigma)\right)= (\E Z_N(\beta))^2 \left(\mathcal{O}\left(\frac 1 {N^2 p^3}\right)+o\left( \frac{1}{Np} \right)\right).
\end{equation}
\end{lemma}
\begin{proof}
We divide the configurations $(\sigma, \tau)$ into ``typical'' pairs of configurations, which by definition,  lie  in the set
\begin{multline}\label{eq:S_N_def}
S_N\coloneqq \{ (\sigma,\tau) \in \{-1,+1\}^N \times \{-1,+1\}^N: \\|\sigma|^2 \le N (Np)^{\frac 15},  |\tau|^2 \le N (Np)^{\frac 15} \mbox{ and }|\sigma \tau|^2 \leq N (Np)^{\frac 1{5}} \}
\end{multline}
 and atypical configurations $(\sigma,\tau) \notin S_N$. Then, by Taylor expanding the exponential in the case of typical spins,
\begin{eqnarray*}
&&\V\left(\sum_\sigma X(\sigma)\right)=\sum_{\sigma, \tau\in \{-1,+1\}^N}\Cov (X(\sigma), X(\tau))\\
&& =\sum_{(\sigma, \tau)\in S_N}\E \eee^{-\beta H(\sigma)}\E \eee^{-\beta H(\tau)}
\left(
\mathcal{O}\left(\frac 1 {N^2 p^3}\right)+\mathcal{O}\Bigl(\bigl(\frac{|\sigma\tau|^2}{N^2}\frac{\beta^2(1-p)}{4p}\bigr)^2\Bigr)+\mathcal{O}\left(\frac {|\sigma|^2+|\tau|^2}{N^3 p^2}\right)\right)\\
&&+\sum_{(\sigma, \tau)\notin S_N} \E \eee^{-\beta H(\sigma)}\E \eee^{-\beta H(\tau)}
\left(\eee^{\frac{|\sigma\tau|^2}{N^2}\frac{\beta^2(1-p)}{4p}+
\mathcal{O}\left(\frac 1 {N^2 p^3} + \frac {|\sigma|^2+|\tau|^2} {N^3p^2}\right)}-1\right.
\\
&& \quad\quad  \left.-\frac{\beta^2 (1-p)}{4N^2 p}  | \sigma \tau|^2+ \mathcal{O}\left(\frac {|\sigma|^2+|\tau|^2}{N^3 p^2}\right)\right).
\end{eqnarray*}

Now on the one hand, for typical spins we have
\begin{eqnarray*}
&&\sum_{(\sigma, \tau)\in S_N}\E \eee^{-\beta H(\sigma)}\E \eee^{-\beta H(\tau)}\left(
\mathcal{O}\left(\frac 1 {N^2 p^3}\right)+ \mathcal{O}\Bigl(\bigl(\frac{|\sigma\tau|^2}{N^2}\frac{\beta^2(1-p)}{4p}\bigr)^2\Bigr)+\mathcal{O}\left(\frac {|\sigma|^2 + |\tau|^2}{N^3 p^2}\right)\right)\\
&\le &
\sum_{(\sigma, \tau)\in S_N}\E \eee^{-\beta H(\sigma)}\E \eee^{-\beta H(\tau)}\left(
\mathcal{O}\left(\frac 1 {N^2 p^3}\right)+ \mathcal{O}\Bigl((Np)^{-2+\frac 25}\Bigr)\right)\\
&\le & \sum_{\sigma, \tau\in \{-1,+1\}^N}\E \eee^{-\beta H(\sigma)}\E \eee^{-\beta H(\tau)}\left(
\mathcal{O}\left(\frac 1 {N^2 p^3}\right)+ \mathcal{O}\Bigl((Np)^{-\frac{8}{5}}\Bigr)\right)\\
&=&(\E Z_N(\beta))^2\left(
\mathcal{O}\left(\frac 1 {N^2 p^3}\right)+ \mathcal{O}\Bigl((Np)^{-\frac{8}{5}}\Bigr)\right).
\end{eqnarray*}
Note that
$ \mathcal{O}\Bigl((Np)^{-\frac{8}{5}}\Bigr)=o\left( \frac{1}{Np} \right)$
for any choice of $p$ such that $Np \to \infty$.

We proceed to the analysis of the atypical pairs of spins. Let $V_N(k,l,m)$ denote the set
$$\{(\sigma,\tau)\in \{-1,+1\}^N\times\{-1,+1\}^N: \; |\sigma|=k,|\tau|= l, \mbox{ and }|\sigma \tau| = m\}
$$
and let $\nu_N(k,l,m) \coloneqq \# V_N(k,l,m)$ be the number of such pairs.
Taking $\sigma=(\sigma_i)_{i=1}^N$ and $\tau=(\tau_i)_{i=1}^N$ independently and uniformly from $\{-1,+1\}^N$, we can regard $(\sigma_i,\tau_i,\sigma_i\tau_i)$, $1\leq i\leq N$, as  i.i.d.\ three-dimensional random vectors with zero mean.
The covariance matrix of these random vectors is the $3\times 3$ identity matrix because
$$
\sigma_i (\sigma_i\tau_i) = \tau_i, \quad \tau_i (\sigma_i\tau_i) = \sigma_i, \quad \sigma_i^2 = \tau_i^2 = (\sigma_i\tau_i)^2 = 1.
$$
By the three-dimensional Local Central Limit Theorem~\cite{Davis1995}, there is a universal constant $C$ such that
\begin{equation}\label{LCLT}
\nu_N(k,l,m) \leq C 2^{2N} N^{-3/2} \eee^{-\frac{k^2}{2N} - \frac{l^2}{2N} - \frac{m^2}{2N}}, \quad (k,l,m)\in\Z^3.
\end{equation}

From here, \eqref{eq:cov_computation1} and  the inequality $1+x\leq e^x$ we get
\begin{align*}
&\left|2^{-2N}\sum_{(\sigma,\tau) \in V_N(k,l,m)}\Cov(X(\sigma), X(\tau))\right|\\
&\leq
2^{-2N}\sum_{(\sigma,\tau) \in V_N(k,l,m)}
\E \eee^{-\beta H(\sigma)}\E \eee^{-\beta H(\tau)}
\left|\eee^{\frac{|\sigma\tau|^2}{N^2}\frac{\beta^2(1-p)}{4p}+
\mathcal{O}\left(\frac 1 {N^2 p^3} + \frac {|\sigma|^2+|\tau|^2} {N^3p^2}\right)}-1\right.
\\
& \quad\quad  \left.-\frac{\beta^2 (1-p)}{4N^2 p}  | \sigma \tau|^2+ \mathcal{O}\left(\frac {|\sigma|^2+|\tau|^2}{N^3 p^2}\right)\right|
\\
&\leq C   N^{-3/2} \eee^{-\frac{k^2}{2N} - \frac{l^2}{2N} - \frac{m^2}{2N}} \eee^{2A_N(\beta) + \frac{\beta}{2} \frac {k^2+l^2} N + \frac 1 {N^2p^2}   \left(2C_{N,1} + C_{N,2}\frac {k^2 + l^2}N\right)}\\
& \qquad \times
\left(\eee^{\frac{m^2}{N^2}\frac{\beta^2(1-p)}{4p}+
\mathcal{O}\left(\frac 1 {N^2 p^3} + \frac {k^2+l^2} {N^3p^2}\right)}+1+\frac{\beta^2 (1-p)}{4N^2 p} m^2+ \mathcal O\left(\frac {k^2+l^2}{N^3 p^2}\right)\right)
\\
&\leq 2C   N^{-3/2} \eee^{-\frac{k^2}{2N} - \frac{l^2}{2N} - \frac{m^2}{2N}} \eee^{2A_N(\beta) + \frac{\beta}{2} \frac {k^2+l^2} N}\cdot
\eee^{\frac{m^2}{N^2}\frac{\beta^2(1-p)}{4p}+
\mathcal{O}\left(\frac 1 {N^2 p^3} + \frac {k^2+l^2} {N^3p^2}\right)}
\\
&\leq 4CN^{-3/2}
\eee^{2A_N(\beta)}
\eee^{-\frac{(1-\beta)(k^2+l^2)}{2N}-\frac{m^2}{2N}} \cdot \eee^{\frac{m^2}{N^2}\frac{\beta^2(1-p)}{4p}+
\mathcal{O}\left(\frac {k^2+l^2} {N^3p^2}\right)},
\end{align*}
for some constant $C>0$ (that may change from line to line) and $N$ sufficiently large. Recalling that $\beta<1$ and using the fact that $Np\to\infty$, we can find a sufficiently small $\delta>0$ such that
$$
\left|\sum_{(\sigma,\tau) \in V_N(k,l,m)}\Cov(X(\sigma), X(\tau))\right|
\leq
4C N^{-3/2} 2^{2N}
\eee^{2A_N(\beta)}
\eee^{-\delta \frac{(k^2+l^2 + m^2)}{2N}},
$$
for all $N$ sufficiently large. Applying this to the pairs $(\sigma,\tau)\in S_N^c$ we obtain
\begin{align*}
\left|\sum_{(\sigma,\tau) \in S_N^c} \Cov(X(\sigma), X(\tau))\right|
\leq
4CN^{-3/2}  2^{2N}
\eee^{2A_N(\beta)} \sum_{\substack{(k,l,m)\in\Z^3\\N^{-1/2}(k,l,m)\in D_N}}
\eee^{-\delta \frac{(k^2+l^2 + m^2)}{2N}},
\end{align*}
where
$$
D_N\coloneqq\{(x,y,z)\in\R^3\colon x^2 > (Np)^{\frac 1 5} \mbox{ or } y^2 > (Np)^{\frac 1 5}  \mbox{ or } z^2 > (Np)^{\frac 1 5}\}.
$$
Estimating the Riemann sum by the Riemann integral and bounding the tail function of the normal distribution, we obtain the crude estimate
$$
N^{-3/2} \sum_{\substack{(k,l,m)\in\Z^3\\N^{-1/2}(k,l,m)\in D_N}}
\eee^{-\delta \frac{(k^2+l^2 + m^2)}{2N}}
\leq C e^{- \delta' (Np)^{1/5}}
$$
for some $\delta'>0$. Taking everything together, we arrive at
$$
\left|\sum_{(\sigma,\tau) \in S_N^c} \Cov(X(\sigma), X(\tau))\right|
\leq
C 2^{2N}
\eee^{2A_N(\beta)} e^{- \delta' (Np)^{1/5}},
$$
Recalling Lemma~\ref{lem:expect}, we finally get the estimate
$$
\left|\sum_{(\sigma,\tau) \in S_N^c} \Cov(X(\sigma), X(\tau))\right|
\leq
C (\E Z_N(\beta))^2 e^{- \delta' (Np)^{1/5}}
=
(\E Z_N(\beta))^2 \mathcal O\left(\frac{1}{N^2p^3}\right).
$$
Putting the estimates for typical and atypical pairs together we obtain the assertion.
\end{proof}

\subsection{Completing the proof of Theorem \ref{CLT_ZN_1}}
We will see that the results obtained so far  suffice for a CLT for
$$
W_N(\beta)\coloneqq \sum_\sigma \E e^{-\beta H(\sigma)}\left[-\beta H(\sigma)- \frac{\beta}{2N} |\sigma|^2\right]
$$
to carry over to a CLT for $Z_N(\beta)-\E Z_N(\beta)$.

We first will show this Central Limit Theorem for $W_N(\beta)$.
Note that by definition of $H$
$$
W_N(\beta)=\sum_\sigma \E e^{-\beta H(\sigma)}\left[\frac \beta {2Np} \sum_{i,j=1}^N \sigma_i \sigma_j \vep_{i,j}- \frac{\beta}{2N} |\sigma|^2\right].
$$
Hence $W_N(\beta)$ is a linear combination of the $\vep_{i,j}$'s and we have $\E W_N(\beta)=0$. However, we need to control the coefficients of the $\vep_{i,j}$'s.  Let us first compute the variance of $W_N(\beta)$. By Lemma~\ref{lem:covB} and Lemma~\ref{lem:lemma_cov} together with the assumption that $Np^2 \to \infty$,
\begin{eqnarray*}
\V(W_N(\beta)) &=& \sum_{\sigma, \tau \in\{\pm1\}^N} \E e^{-\beta H(\sigma)}\E e^{-\beta H(\tau)}
\frac{\beta^2 (1-p)}{4N^2 p}  | \sigma \tau|^2\\
&=&
\eee^{\frac{\beta^2 (1-p)}{4p}+o(1)}
\sum_{\sigma, \tau \in\{\pm1\}^N}
\frac{\beta^2 (1-p)}{4N^2 p}  | \sigma \tau|^2
\eee^{\frac{\beta}{2} \frac {|\sigma|^2+|\tau|^2} N + \frac 1 {N^2p^2}   \left(C_{N,2}\frac {|\sigma|^2 + |\tau|^2}N\right)}\\
&=&
\eee^{\frac{\beta^2 (1-p)}{4p}+o(1)}\frac{\beta^2 (1-p)}{4N p}
\sum_{\sigma, \tau \in\{\pm1\}^N}
\frac{1}{N}  | \sigma \tau|^2
\eee^{\frac{\beta}{2} \frac {|\sigma|^2+|\tau|^2} N}.
\end{eqnarray*}

Now again, under the uniform distribution on $\{\pm1\}^N\times \{\pm1\}^N$ the random vector $(\frac{|\sigma|}{\sqrt N},\frac{|\tau|}{\sqrt N},\frac{|\sigma \tau|}{\sqrt N})$ converges in distribution to a three-dimensional Gaussian random vector $(\xi_1,\xi_2,\xi_3)$ with identity covariance matrix. Moreover, since $\beta <1$ the functions
\begin{equation} \label{eq:UI}
G_N(\sigma, \tau)\coloneqq \frac{1}{N}  | \sigma \tau|^2
\eee^{\frac{\beta}{2} \frac {|\sigma|^2+|\tau|^2} N}
\end{equation}
are uniformly integrable with respect to this uniform distribution on $\{\pm1\}^N\times \{\pm1\}^N$. This can exactly be shown as in \cite[Proof of Theorem V.9.4]{Ellis-EntropyLargeDeviationsAndStatisticalMechanics}. Note that as $\beta <1$ and by the independence  of $\xi_1, \xi_2$ and $\xi_3$ we have
$$
 \E_{\xi_1,\xi_2,\xi_3}\left(\xi_3^2 \exp\left(\frac{\beta}2(\xi_1^2+\xi_2^2)\right)\right) = \left(\mathbb E_{\xi_1} \left( e^{\frac{\beta}{2} \xi_1^2}\right)\right)^2=\frac{1}{1-\beta}.
$$
Thus
\begin{eqnarray}\label{eq:VWN}
\V(W_N(\beta)) &= & 2^{2N} \eee^{\frac{\beta^2 (1-p)}{4p}+o(1)}\frac{\beta^2 (1-p)}{4N p} \frac 1{2^{2N}}
\sum_{\sigma, \tau \in\{\pm1\}^N}
\frac{1}{N}  | \sigma \tau|^2
\eee^{\frac{\beta}{2} \frac {|\sigma|^2+|\tau|^2} N} \nonumber\\
&\sim &   \frac{\beta^2 (1-p)}{4N p}  e^{\frac{\beta^2 (1-p)}{4p}} 4^N \E_{\xi_1,\xi_2,\xi_3}\left(\xi_3^2 \exp\left(\frac{\beta}2(\xi_1^2+\xi_2^2)\right)\right)
\nonumber\\
&\sim &  \frac{1}{1-\beta}  \frac{\beta^2 (1-p)}{4N p}  e^{\frac{\beta^2 (1-p)}{4p}} 4^N\nonumber\\
&\sim&    \frac{\beta^2 (1-p)}{4N p} (\E Z_N(\beta))^2
\end{eqnarray}
  according to Lemma \ref{lem:expect}.

To prove the CLT for $W_N(\beta)$ let us define
$$
X_{i,j} \coloneqq  \vep_{i,j} \sum_{\sigma} \frac{\beta}{2Np}\sigma_i \sigma_j \E \eee^{-\beta H(\sigma)}
$$
Then
\begin{eqnarray*}
W_N(\beta) &=& \sum_\sigma \E e^{-\beta H(\sigma)}\left[-\beta H(\sigma)- \frac{\beta}{2N} |\sigma|^2\right]\\
&=:& 
\sum_{i,j} (X_{i,j}-\E X_{i,j}).
\end{eqnarray*}
We shall verify the conditions of the Lyapunov CLT. Note that by  Lemma \ref{lem:expect} we have that $\E Z_N (\beta)\sim \eee^{\frac{(1-p)\beta^2}{8p}} \frac {2^N} {\sqrt{1-\beta}}$ and therefore
$$
 |X_{i,j}|  \leq  \Big| \vep_{i,j} \sum_{\sigma} \frac{\beta}{2Np}\sigma_i \sigma_j \E \eee^{-\beta H_N(\sigma)} \Big| \leq C 2^N \frac{1}{Np}  \eee^{\frac{\beta^2 (1-p)}{8  p}}
$$
for yet another constant $C>0$, that depends on $\beta$.
This means, recalling that $(1-p) Np \to \infty$,
\begin{equation}\label{est_X_ij}
|X_{i,j}| =
o\left(\frac{ \sqrt{(1-p)}}{  \sqrt{N p}}  e^{\frac{\beta^2 (1-p)}{8p}} 2^N \right) \quad \text{ resp. }\quad |X_{i,j}| =
o(\sqrt{\V(W_N(\beta))}),
\end{equation}
where we used $\V(W_N(\beta)) \sim C' \frac{\beta^2 (1-p)}{4N p}  e^{\frac{\beta^2 (1-p)}{4p}} 4^N$.
This means that Lyapunov's condition is satisfied. 
Indeed, take any $\delta >0$ and write $s_N\coloneqq\sqrt{\V(W_N(\beta))}$. Then by \eqref{est_X_ij} we have  $|X_{i,j}-\mathbb E(|X_{i,j}|)|^{\delta}=o(s_N^{\delta})$ and hence
\begin{equation*}
\frac 1{s_N^{2+\delta}}\sum_{i,j=1}^N \mathbb E[|X_{i,j}-\E X_{i,j}|^{2+\delta}]
= o(1) \frac{1}{s_N^2}\sum_{i,j=1}^N \mathbb E[|X_{i,j}-\E X_{i,j}|^{2}]\\
= o(1).
\end{equation*}

Therefore $W_N(\beta)$ satisfies the CLT. Taking into account that $\E W_N(\beta)=0$ and the asymptotic variance of $W_N(\beta)$ computed in \eqref{eq:VWN} this means
$$
\frac 2 \beta \sqrt{\frac{pN}{1-p}} \frac{W_N(\beta)}{\E Z_N(\beta)} \to \mathcal{N}(0,1)
$$
in distribution.  The final observation is that due to what we computed before in \eqref{estimate_var_sum_X}
$$
\left|\left|\sqrt{\frac{pN}{1-p}} \frac{\sum_\sigma X(\sigma)}{\E Z_N(\beta)}\right|\right|_2^2 =\mathcal{O}\left(\frac 1 {Np^2}\right) +o(1) = o(1)
$$
by the assumption  $p^2 N \to \infty$ made in the statement of the theorem.
Therefore, recalling that $\sum_{\sigma} X(\sigma) = Z_N(\beta)- \mathbb E(Z_N(\beta))-W_N(\beta)$, the Central Limit Theorem for $W_N(\beta)$ carries over to $Z_N(\beta)-\E Z_N(\beta)$ and we have
$$
\sqrt {\frac{pN}{1-p}} \left( \frac{Z_N(\beta)}{\E Z_N(\beta)} -1\right) \to \mathcal{N}\left(0,\frac {\beta^2} 4\right).
$$
in distribution, which is the assertion.
\hfill $\Box$

\section{Proof of Theorem \ref{CLT_ZN_2}}
\subsection{Method of proof}
Very roughly, the idea of our proof of Theorem~\ref{CLT_ZN_2} can be explained as follows. The $2^N$ random variables $H_N(\sigma)$, $\sigma\in \{-1,+1\}^N$, are clearly dependent.
We will  look for a representation of the form
$$
H_N(\sigma) = \hat H_N(\sigma) + \Delta_N,
$$
where $\Delta_N$ is certain random variable containing ``most'' of the dependence between the $H_N(\sigma)$'s, and the random variables $\hat H_N(\sigma)$, $\sigma\in \{-1,+1\}^N$, behave as ``almost'' independent. The partition function can be represented as
$$
Z_N(\beta) = \sum_{\sigma \in \{-1,+1\}^N}\exp(-\beta H_N(\sigma)) = e^{-\beta \Delta_n}\sum_{\sigma \in \{-1,+1\}^N}\exp(-\beta \hat H_N(\sigma)).
$$
We will show that the fluctuations of the sum on the right-hand side around its mean are negligible (that is, the sum is ``almost'' deterministic) and most of the fluctuations of $Z_N(\beta)$ come from the term $e^{-\beta \Delta_N}$.
To make this approach work, it is necessary to identify $\Delta_N$. In a search for $\Delta_N$, it is natural to try out simple symmetric functions of the random variables $\eps_{ij}$, the most simple examples being the following ones:
$$
\xi_N \coloneqq \frac 1 {\sqrt{N p(1-p)}} \sum_{i=1}^N \vep_{i,i}
\;\;
\text{ and }
\;\;
\eta_N \coloneqq \frac 1 {N \sqrt{p(1-p)}} \sum_{i,j=1}^N \vep_{i,j}.
$$
Making the Ansatz $\Delta_N:= \alpha_N \xi_N + \beta_N \eta_N$ with unknown deterministic coefficients $\alpha_N$ and $\beta_N$ we have found, after some trials and errors, that the ``correct'' form of $\Delta_N$ is
$$
\Delta_N = -\frac{\sqrt{1-p}}{2\sqrt {Np}}\xi_N - \frac{\beta (1-2p)\sqrt{1-p}}{8N p^{3/2}}\eta_N.
$$

It is now time to introduce some notation. Let $\gamma\coloneqq \frac{\beta}{2 Np}$. For each $\sigma\in \{-1,+1\}^N$ define the random variable
$$
T(\sigma)\coloneqq \exp\left(\gamma \left( \sum_{i,j=1}^N \vep_{i,j}\sigma_i \sigma_j\right) - \frac{\beta\sqrt{1-p}}{2\sqrt {Np}}\xi_N-
\frac{\beta^2 (1-2p)\sqrt{1-p}}{8N p^{3/2}}\eta_N	\right).
$$
Define the corresponding modified partition function by
$$
\hat Z_N(\beta) := \sum_{\sigma\in \{-1,+1\}^N} T(\sigma)
= Z_N(\beta) \exp\left(- \frac{\beta\sqrt{1-p}}{2\sqrt {Np}}\xi_N-
\frac{\beta^2 (1-2p)\sqrt{1-p}}{8N p^{3/2}}\eta_N	\right).
$$
Using this notation,  we can write
\begin{equation}\label{eq:Z_N_hat_Z_N}
Z_N(\beta)
=
\hat Z_N(\beta) \exp\left(\frac{\beta\sqrt{1-p}}{2\sqrt {Np}}\xi_N
+\frac{\beta^2 (1-2p)\sqrt{1-p}}{8N p^{3/2}}\eta_N	\right).
\end{equation}
The key result of our proof is Proposition~\ref{var_sum} which states that
$$
\V \left(\frac{\hat Z_N(\beta)}{\E \hat Z_N(\beta)}\right) = o\left(\frac 1 {Np}\right).
$$
This means that $\hat Z_N(\beta)$ is an almost deterministic quantity. The fluctuations of $Z_N(\beta)$ are therefore determined by the fluctuations of $e^{-\beta \Delta_N}$, which is the second factor on the right-hand side of~\eqref{eq:Z_N_hat_Z_N}. Note that by the central limit theorem,  both $\xi_N$ and $\eta_N$ converge to a standard normal random variable after centering, so that the asymptotic  fluctuations of $\Delta_N$ are Gaussian with an easily identifiable variance. It turns out that in the regime when $p^3 N^2 \to \infty$, the variance converges to $0$, which means that we can use the Taylor expansion  $e^{-\beta \Delta_N} = 1 - \beta \Delta_N + o(\Delta_N)$ together with the central limit theorem for $\xi_N$ and $\eta_N$ to  prove Theorem~\ref{CLT_ZN_2}. Let us mention that this approach allows us to identify the quantity which causes the fluctuations of the partition function: it is either $\eta_N$ (in the regime when $p^2N  \to 0$ but $p^3 N^2 \to \infty$) or $\xi_N$ (in the regime when $p^2 N\to\infty$). The latter regime will not be considered here because it has been already treated in Section~\ref{sec:proof_CLT_ZN_1} by a different method. Finally, in the boundary case when $p^2 N  \to \text{const}$, both terms $\xi_N$ and $\eta_N$ contribute to the fluctuations of $Z_N(\beta)$. 

We believe that the applicability of our method is not restricted to the model studied here. Some other examples with a similar behaviour can be found in~\cite{kabluchko_levels}.


\subsection{Expectations and covariances}
To make the above considerations rigorous, we have to provide asymptotic expressions for $\E \hat Z_N(\beta)$ and $\V \hat Z_N(\beta)$. This will be done in a sequence of lemmas.   In the following, the inverse temperature is fixed and satisfies $0<\beta <1$.

\begin{lemma}\label{lem:expect1}
For all $\sigma \in \{-1,+1\}^N$ we have
$$
\E T(\sigma)= \exp\left(N \log(1-p+p \eee^{-\frac{\gamma^2}2 (1-2p)})+ \sum_{\substack{i,j=1\\i\neq j}}^N \log(1-p+p\eee^{\gamma \sigma_i \sigma_j-\frac{\gamma^2}2 (1-2p)})\right).
$$
\end{lemma}
\begin{proof}
This is an elementary computation using
\begin{equation*}
T(\sigma)= \exp \left( -\frac{\gamma^2}{2}(1-2p) \sum_{i=1}^N \varepsilon_{i,i} + \sum_{\substack{i,j=1\\i\neq j}}^N \left( \gamma \sigma_i \sigma_j - \frac{\gamma^2}{2}(1-2p) \right) \varepsilon_{i,j}\right)
\end{equation*}
together with the independence of the $(\vep_{i,j})_{i,j}$ and the definitions of $\gamma$, $\xi_N$, and $\eta_N$.
\end{proof}
From here we prove
\begin{lemma}\label{lem:expect2}
For all $\sigma \in \{-1,+1\}^N$, $\E T(\sigma)$ has the following asymptotic behavior when $p \to 0$ and $N^2p^3 \to \infty$:
$$
\E T(\sigma)= \exp\left(\frac{\beta^2}8\left(1-\frac 1 {Np}\right)-\frac{\beta^4}{12\cdot 16 N^2p^3}+\frac{\beta}2\left(\frac{|\sigma|^2}N -1\right)+
\mathcal{O}(\gamma^3 p |\sigma|^2)+o\left(\frac 1 {Np}\right)\right).
$$
In this and similar results, the rest terms only depend on $\sigma$ (if at all) in the way specified in the corresponding brackets.
\end{lemma}

\begin{proof}
Let us define
$$
f_1\coloneqq \log(1-p+p \eee^{-\frac{\gamma^2}2 (1-2p)})
$$
and
$$
f_2(x)\coloneqq\log(1-p+p\eee^{\gamma x-\frac{\gamma^2}2 (1-2p)}).
$$
Then by Lemma \ref{lem:expect1} we have
\begin{equation*}
\E T(\sigma)= \exp\Bigg(N f_1 + \sum_{\substack{i,j=1\\i\neq j}}^N f_2(\sigma_i \sigma_j)\Bigg).
\end{equation*}
Expanding the exponential and the logarithm in the definition of $f_1$, and taking into account that we are in a regime where
$p \to 0$, $\gamma \to 0$, and $Np \to \infty$, we obtain
$$
N f_1 = -Np\frac{\gamma^2}2 +o\left(\frac{1}{Np}\right)= -\frac{\beta^2}{8Np}+o\left(\frac{1}{Np}\right).
$$
On the other hand let us consider
$\sum_{i \neq j}^N f_2(\sigma_i\sigma_j).$
Recall that the product $\sigma_i\sigma_j$ can only take the values $\pm 1$.  At these values we may rewrite $f_2$ as
\begin{equation}\label{expansion_f_2}
f_2(\sigma_i\sigma_j)=a_0+a_1 \sigma_i \sigma_j
\end{equation} leading to
$$
\sum_{\substack{i,j=1\\i\neq j}}^N f_2(\sigma_i\sigma_j)= N(N-1) a_0+ a_1(|\sigma|^2-N),
$$
where again we write $|\sigma|^2=(\sum_{i=1}^N \sigma_i)^2$.
 Solving \eqref{expansion_f_2} for $a_0$ and $a_1$ gives
 \begin{equation*}a_0=\frac {f_2(1)+f_2(-1)}{2} \quad \text{and} \quad a_1=\frac {f_2(1)-f_2(-1)}{2}.
 \end{equation*}
Taking into account that, with the notation given by \eqref{eq:centralfunc}, we have $f_2(\pm 1)=F_1(\pm 1,p,\gamma)$
 and using the expansion
of $F_1$ given in \eqref{f1a} we obtain:
$$
a_0=  \frac{\gamma^2 p^2}2-\frac 1{12}p\gamma^4 +\mathcal{O}(p^2 \gamma^4)+\mathcal{O}(\gamma^6 p).
$$
Hence
$$
N(N-1) a_0= \frac{\beta^2}8-\frac{\beta^4}{12\cdot 16 N^2p^3}+o\left(\frac 1{Np}\right),
$$
where we have used that $p \to 0 $ and $N^2p^3 \to \infty$.

Similarly, using~\eqref{f1b} we compute
$$
a_1=  \gamma p + \mathcal{O}(p\gamma^3).
$$
Now,
$$
\gamma p \, (|\sigma|^2-N)=\frac{\beta}{2}\left(\frac{|\sigma|^2}N-1\right).
$$
From this observation and the above expansion of $a_0$ the assertion follows.
\end{proof}

In a similar fashion we can now compute the covariances.
\begin{lemma}\label{lem:cov1}
For all $\sigma, \tau \in \{-1,+1\}^N$ we have
$$
\E[T(\sigma) T(\tau)]=
\exp\left(N \log(1-p+p \eee^{-\gamma^2 (1-2p)})+ \sum_{\substack{i,j=1\\i\neq j}}^N \log(1-p+p\eee^{\gamma(\sigma_i \sigma_j+\tau_i\tau_j)-\gamma^2 (1-2p)})\right).
$$
\end{lemma}
\begin{proof}
This is almost the same computation as in the proof of Lemma \ref{lem:expect1}.
\end{proof}
\begin{lemma}\label{lem:cov2}
For all $\sigma, \tau \in \{-1,+1\}^N$ and for $p \to 0$, $p^3N^2 \to \infty$ we have
\begin{multline*}
\E[T(\sigma) T(\tau)]=
\exp\left(\frac{\beta^2}4\left(1-\frac{1}{Np}\right)-\frac{\beta^4}{6\cdot 16 N^2p^3}+\frac{\beta}2\left(\frac{|\sigma|^2}N+\frac{|\tau|^2}N -2\right)
\right.
\\
\left.
+
\mathcal{O}(\gamma^3 p (|\sigma|^2+|\tau|^2))+ \frac{\beta^2(1-p)}{4Np}\left(\frac{|\sigma \tau|^2}N-1\right)+\mathcal{O}(\gamma^4 p|\sigma \tau|^2)+
o\left(\frac 1 {Np}\right)\right).
\end{multline*}
\end{lemma}
\begin{proof}
The proof is similar to the proof of Lemma \ref{lem:expect2}.
Starting from Lemma \ref{lem:cov1} we set
$$
f_3\coloneqq \log(1-p+p \eee^{-\gamma^2 (1-2p)})
$$
and
$$
f_4(\sigma_i \sigma_j+\tau_i\tau_j)\coloneqq \log(1-p+p\eee^{\gamma (\sigma_i \sigma_j+\tau_i\tau_j)-\gamma^2 (1-2p)}),
$$
which implies
\begin{equation*}
\E[T(\sigma) T(\tau)]= \exp \left( N f_3 + \sum_{\substack{i,j=1\\i\neq j}}^{N} f_4(\sigma_i \sigma_j+\tau_i\tau_j) \right).
\end{equation*}
An easy expansion as in the proof of Lemma \ref{lem:expect2} shows that
$$N f_3= -\frac{\beta^2}{4Np}+o\left(\frac 1{Np}\right).$$

Again similar to the proof of Lemma \ref{lem:expect2} we observe that $\sigma_i \sigma_j+\tau_i\tau_j \in \{-2,0,2\}$ for all possible choices of $\sigma$ and $\tau$.
For these values we represent $f_4$ as
$$
f_4(x+y)= b_0+b_1 x +b_2 y +b_{12}xy \qquad x,y \in \{-1,+1\}
$$
with coefficients $b_0, b_1,b_2,b_{12}$ to be determined now.
These are readily computed to be given by
$$
b_0= \frac{f_4(2)+f_4(-2)+2f(0)}4, \quad b_1=b_2=\frac{f_4(2)-f_4(-2)}4,
$$
and
$$
b_{12}= \frac{f_4(2)+f_4(-2)-2f(0)}4.
$$
Note that with the notation introduced in \eqref{eq:centralfunc} we have $f_4(\pm 2)=F_2(\pm 2, p, \gamma)$ as well as
 $f_4(0)= F_2(0,p,\gamma)$. Therefore, using the expansion \eqref{f2a} we obtain the following expansions for the coefficients
$$
N(N-1) b_0= \frac{\beta^2}4-\frac{\beta^4}{6\cdot 16 N^2p^3}+o\left(\frac 1{Np}\right).
$$
In particular,  up to $o\left(\frac 1 {Np}\right)$-terms we have $N(N-1)b_0=2N(N-1)a_0$ with $a_0$ as in  the proof of Lemma \ref{lem:expect2}.

Moreover, by~\eqref{f2b},
$$b_{12}= \gamma^2 p(1-p)-\frac 23 \gamma^4 p +\mathcal{O}(\gamma^4 p^2) + \mathcal{O}(\gamma^6 p).$$ Taking into account that $|\sigma \tau|^2 \le N^2$, $N^2 p^3 \to \infty$, and the definition of $\gamma$, this implies
$$
b_{12}(|\sigma\tau|^2-N)= \frac{\beta^2(1-p)}{4Np}\left(\frac{|\sigma \tau|^2}N-1\right)-\frac 23 \gamma^4 p|\sigma \tau|^2+o\left(\frac 1{Np}\right).
$$

Finally, again similar to the proof of Lemma \ref{lem:expect2}
$$
b_1= b_2=\gamma p+ \mathcal O( p\gamma^3).
$$
Putting these estimates together yields the assertion.
\end{proof}

To make use of these results we will need to compute the expectation of $\sum_\sigma T(\sigma)$. To this end recall the set of typical pairs of spins from Equation~\eqref{eq:S_N_def}
$$
S_N\coloneqq\{(\sigma,\tau)\in \{\pm 1\}^N \times \{\pm 1\}^N: |\sigma|^2 \le N (Np)^{\frac 15},|\tau|^2 \le N (Np)^{\frac 15},
|\sigma\tau|^2\le N (Np)^{\frac{1}{5}}\}.
$$
In a slight abuse of notation we will also say that a configuration $\sigma$ is typical and write $\sigma \in S_N$, when
$ |\sigma|^2 \le N (Np)^{\frac 15}$.

\begin{lemma} \label{lem:expect_sum}
In the regime when $p \to 0$ and $N^2p^3 \to \infty$ we have the following asymptotic behaviour for the expectation of $\hat Z_N(\beta) = \sum_\sigma T(\sigma)$:
$$
\E \hat Z_N(\beta)
=
2^N \eee^{\frac{\beta^2}{8}-\frac{\beta}{2}}\frac{1}{\sqrt{1-\beta}} \left(1+o(1) \right).
$$
\end{lemma}
\begin{rem}
Note that in particular we have $C_1 2^N \leq \E \hat Z_N(\beta)\leq C_2 2^N$ for suitable constants $0<C_1<C_2<\infty$.
\end{rem}
\begin{proof}[Proof of Lemma \ref{lem:expect_sum}]
According to Lemma \ref{lem:expect2} we have
\begin{align}\label{Esplit}
&\E \sum_\sigma T(\sigma) = \E \sum_{\sigma \in S_N} T(\sigma)+\E \sum_{\sigma \notin S_N} T(\sigma)\nonumber\\
&= \sum_{\sigma \in S_N} \exp\left(\frac{\beta^2}8\left(1-\frac 1 {Np}\right)-\frac{\beta^4}{12\cdot 16 N^2p^3}+\frac{\beta}2\left(\frac{|\sigma|^2}N -1\right)+o\left(\frac 1 {Np}\right)\right)\nonumber\\
&+ \sum_{\sigma \notin S_N} \exp\left(\frac{\beta^2}8\left(1-\frac 1 {Np}\right)-\frac{\beta^4}{12\cdot 16 N^2p^3}+\frac{\beta}2\left(\frac{|\sigma|^2}N -1\right)+
\mathcal{O}(\gamma^3 p |\sigma|^2)+o\left(\frac 1 {Np}\right)\right),
\end{align}
where we used for the sum over  $\sigma \in S_N$ that the error term
$\mathcal{O}(\gamma^3 p |\sigma|^2)$ can be absorbed in the $o(\frac 1 {Np})$-term.
We start with the first term, summing over all $\sigma \in S_N$.
Note that we have
\begin{equation*}
\eee^{\frac{-\beta^2}{8Np}+\frac{-\beta^4}{12\cdot 16 N^2p^3}+o\left(\frac 1 {Np}\right)}
=
1+\mathcal{O}\left(\frac 1 {Np}\right)+\mathcal{O}\left(\frac 1{N^2p^3}\right)+o\left(\frac 1 {Np}\right)=1+o(1),
\end{equation*}
since $Np\to \infty$ and $N^2p^3 \to \infty$ by assumption.
Hence we have
\begin{equation*}
\E \sum_{\sigma \in S_N} T(\sigma)=\left( 1+o(1) \right)\eee^{\frac{\beta^2}8 -\frac{\beta} 2}\sum_{\sigma \in S_N}\eee^{\frac{\beta}2 \frac{|\sigma|^2}N} .
\end{equation*}
Recall that under the uniform measure on $\{-1,+1\}^N$ the random variable $|\sigma|/\sqrt N$ converges in distribution to a standard normal random variable. Due to uniform integrability of the sequence of random variables $\sigma \mapsto e^{\beta \frac{|\sigma|^2}{2N}}$, $N\in\N$,  for $\beta<1$ (see~\cite[Proof of Theorem V.9.4]{Ellis-EntropyLargeDeviationsAndStatisticalMechanics}) this implies that $2^{-N}\sum_{\sigma \in \{-1,1\}^N}\eee^{\frac{\beta}2 \frac{|\sigma|^2}N}$
  is asymptotic to $\frac 1 {\sqrt{2\pi}} \int_{-\infty}^\infty
\eee^{-(1-\beta)\frac{x^2}2}dx= \frac 1 {\sqrt{1-\beta}}$. Indeed, then also $2^{-N}\sum_{\sigma \in S_N}\eee^{\frac{\beta}2 \frac{|\sigma|^2}N}$ is asymptotic to $ \frac 1 {\sqrt{1-\beta}}$, since we shall see in~\eqref{rest_int} below that $2^{-N}\sum_{\sigma \notin S_N}\eee^{\frac{\beta}2 \frac{|\sigma|^2}N}=o(1)$.
Hence,
\begin{align*}
\E \sum_{\sigma \in S_N} T(\sigma)=2^N\eee^{\frac{\beta^2}8 -\frac{\beta} 2}\frac{1}{\sqrt{1-\beta}}\left( 1+ o(1) \right).
\end{align*}

So we just need to show that the second summand in \eqref{Esplit} over the atypical $\sigma$ is $o(2^N)$. To this end observe that
\begin{align*}
&\sum_{\sigma \notin S_N} \exp\left(\frac{\beta^2}8\left(1-\frac 1 {Np}\right)-\frac{\beta^4}{12\cdot 16 N^2p^3}+\frac{\beta}2\left(\frac{|\sigma|^2}N -1\right)+
\mathcal{O}(\gamma^3 p |\sigma|^2)+o\left(\frac 1 {Np}\right)\right) \\
\le& K_1 \sum_{\sigma \notin S_N} \eee^{\frac{\beta}2\frac{|\sigma|^2}N+
K_2\gamma^3 p |\sigma|^2}= K_1 \sum_{k: k^2 > N (Np)^{\frac 15}}\sum_{\sigma: |\sigma|^2=k^2} \eee^{\frac{\beta}2\frac{k^2}N+
K_2\gamma^3 p k^2}
\end{align*}
for constants $K_1,K_2>0$.
By the Stirling formula, there is constant $K_3>0$ such that
$$\# \{\sigma:|\sigma|^2=k^2\}\le K_3\frac{2^N}{\sqrt N} \eee^{\frac{-k^2}{2N}}.$$
Plugging this into the above estimate we arrive at
\begin{eqnarray*}
K_1 \sum_{k: k^2 > N (Np)^{\frac 15}}\sum_{\sigma: |\sigma|^2=k^2} \eee^{\frac{\beta}2\frac{k^2}N+
K_2\gamma^3 p k^2}
&\le& K_4 \sum_{k: k^2 > N (Np)^{\frac 15}} \frac{2^N}{\sqrt N} \eee^{\frac{-k^2}{2N}+\frac{\beta}2\frac{k^2}N+
K_4\frac1 {N^2p^2} \frac{k^2} N}\\
&\le & \sum_{k: k^2 > N (Np)^{\frac 15}} \frac{2^N}{\sqrt N} \eee^{-\delta \frac{k^2}{2N}}
\end{eqnarray*}
for some appropriate $K_4>0$ and $\delta>0$, whenever $N$ is large enough. By comparing the sum to the Riemann integral
$2^N \int_{\frac 12 (Np)^{\frac 15}}^\infty \eee^{-\delta \frac{x^2}2}dx$ we see that
\begin{equation}\label{rest_int}
  \sum_{\sigma \notin S_N} \eee^{\frac{\beta}2\frac{|\sigma|^2}N+
K_2\gamma^3 p |\sigma|^2} =2^N o(1)
\end{equation}
and hence
\begin{align*}
& \sum_{\sigma \notin S_N} \exp\left(\frac{\beta^2}8\left(1-\frac 1 {Np}\right)-\frac{\beta^4}{12\cdot 16 N^2p^3}+\frac{\beta}2\left(\frac{|\sigma|^2}N -1\right)+
\mathcal{O}(\gamma^3 p |\sigma|^2)+o\left(\frac 1 {Np}\right)\right)
\\
&=o(2^N).
\end{align*}
 \end{proof}

\subsection{Estimate for the variance}
Next we can show that $\hat Z_N(\beta)$ has a small variance compared to $(\E \hat Z_N(\beta))^2$.
\begin{proposition}\label{var_sum}
In the regime when $p \to 0$ and $p^3N^2 \to \infty$, the variance of $\hat Z_N(\beta) = \sum_\sigma T(\sigma)$ satisfies
$$
Np \cdot \V \hat Z_N(\beta) = o\left(4^N\right).
$$
Consequently,
$$
\V \left(\frac {\hat Z_N(\beta)}{\E \hat Z_N(\beta)}\right) = o\left(\frac 1 {Np}\right).
$$
\end{proposition}

\begin{proof}
With the above definition of the set of typical pairs of spins $S_N$ note that
\begin{eqnarray*}
\V \sum_\sigma T(\sigma) &=& \sum_{\sigma, \tau} \Cov(T(\sigma), T(\tau))\\
&=& \sum_{(\sigma, \tau)\in S_N} \Cov(T(\sigma), T(\tau))+\sum_{(\sigma, \tau) \notin S_N} \Cov(T(\sigma), T(\tau)).
\end{eqnarray*}
Let us treat the second summand on the right-hand side first.
To this end, denote by $V_N(k,l,m)$ the set of pairs
$$(\sigma,\tau)\in \{-1,+1\}^N\times\{-1,+1\}^N, \quad \mbox{for which }|\sigma|=k,|\tau|= l, \mbox{and }|\sigma \tau| = m.$$
Moreover, by $\nu_N(k,l,m) := \# V_N(k,l,m)$ let us denote the number of such pairs.
We again want to apply a Local Limit Theorem.
To this end
recall that from \eqref{LCLT} we have that
$$
\nu_N(k,l,m) \leq C 2^{2N} N^{-3/2} \eee^{-\frac{k^2}{2N} - \frac{l^2}{2N} - \frac{m^2}{2N}}, \quad (k,l,m)\in\Z^3.
$$
First we estimate the contribution of the untypical pairs of spins. 
Using Lemma \ref{lem:cov2} we see that for any $(k,l,m)\in\Z^3$ there are universal constants $C_1, C_2 >0$ such that
\begin{multline*}
2^{-2N}\sum_{(\sigma,\tau) \in V_N(k,l,m)}
\E [T(\sigma) T(\tau)]
\leq \\
C_1N^{-3/2}
\exp\left(\frac{\beta-1}2\left(\frac{k^2}N+\frac{l^2}N\right)
+
C_2\left(\frac{\beta^3}{N^2p^2}\frac{k^2+l^2}N\right)+\left(\frac{\beta^2}{4 Np}-\frac{1}{2}\right)\frac{m^2}N\right),
\end{multline*}
where we used that $Np\to\infty$, $N^2p^3\to\infty$ and $|m|\leq N$. 
Since $\beta<1$ and $Np\to\infty$,  we obtain for some $\delta >0$ and $N$ large enough
\begin{equation*}
2^{-2N}\sum_{(\sigma,\tau) \in V_N(k,l,m)}
\E [T(\sigma) T(\tau)]
\leq
C_1 N^{-3/2}
\eee^{-\delta \frac{k^2+l^2 + m^2} {2N}}.
\end{equation*}

In a very similar way we can use Lemma~\ref{lem:cov1} to bound the sum of the $\E [T(\sigma)]\E[T(\tau)]$-terms as follows:
\begin{align*}
&2^{-2N}\sum_{(\sigma,\tau) \in V_N(k,l,m)}\E [T(\sigma)]\, \E[T(\tau)]
\le  C_1 N^{-3/2}
\eee^{-\delta \frac{k^2+l^2 + m^2} {2N}}
\end{align*}
for some $\delta >0$ and $N$ sufficiently large (and possibly a different constant $C_1$).

We can therefore conclude that there is a constant $C_1>0$ such that
$$
2^{-2N}\sum_{(\sigma,\tau) \in V_N(k,l,m)}
\left| \Cov (T(\sigma), T(\tau))\right|
\leq
C_1 N^{-3/2} \eee^{-\delta \frac{k^2+l^2 + m^2} {2N}}
$$
for some $\delta >0$ and $N$ sufficiently large.
Hence, the contribution of the untypical spins can be bounded above as follows:
\begin{equation*}
\left|\sum_{(\sigma,\tau) \in S_N^c}
 \Cov (T(\sigma), T(\tau)) \right|
\leq
C_1 4^{N} N^{-3/2}
\sum_{\substack{(k,l,m)\in\Z^3\\N^{-1/2}(k,l,m)\in D_N}} \eee^{-\delta \frac{k^2+l^2 + m^2} {2N}},
\end{equation*}
where $D_N\coloneqq\{(x,y,z)\in\R^3\colon |x|>(Np)^{\frac 1{10}} \text{ or } |y|>(Np)^{\frac 1{10}} \text{ or } |z|>(Np)^{\frac 1{10}} \}$.
Next we demonstrate that the right hand is $o(4^N / (Np))$.
We consider the sum on the right-hand side as a Riemann sum. We will bound it from above by the corresponding
integral over a larger domain. Including the pre-factor $N^{-3/2}$ this apporoach yields
$$
 N^{-3/2}
\sum_{\substack{(k,l,m)\in\Z^3\\N^{-1/2}(k,l,m)\in D_N}} \eee^{-\delta \frac{k^2+l^2 + m^2} {2N}}
\leq
\int_{\frac 13 D_N} \eee^{-\delta \frac{x^2+y^2 + z^2} {2}} \dd x\, \dd y \,  \dd z .
$$
We use the classical bound
$$
\int_{\lambda}^\infty\eee^{-\delta \frac{x^2} {2}} \dd x \le \frac {C_3} {\delta \lambda}
e^{-\delta \frac{\lambda^2}2}.
$$
together with the fact  that on $D_N$ at least one of $x^2,y^2$, and $z^2$  has to be larger than $(Np)^{\frac 1{5}}$. Hence we see that the above integral satisfies
$$
\int_{\frac 13 D_N} \eee^{-\delta \frac{x^2+y^2 + z^2} {2}} \dd x\, \dd y \,  \dd z \le
C_4 \frac 1 {(Np)^{\frac 1{10}}}e^{-C_5 (Np)^{\frac 15}}
$$
for some constants $C_4>0$ and $C_5>0$.
Since $Np \to \infty$ as $N\to\infty$ this implies that
$$
\int_{\frac 13 D_N} \eee^{-\delta \frac{x^2+y^2 + z^2} {2}} \dd x\, \dd y \,  \dd z=o\left(\frac 1 {Np}\right).
$$
Altogether, this shows that
\begin{equation}\label{est_atyp_cov}
Np \sum_{(\sigma, \tau) \in S_N^c}  \Cov (T(\sigma), T(\tau))
=
o(4^N).
\end{equation}
We are now going to prove the analogue of \eqref{est_atyp_cov} for the sum over the typical pairs of spins $(\sigma, \tau) \in S_N$.
Again, we use Lemmas \ref{lem:cov2} and \ref{lem:expect2} and we observe that for
$(\sigma,\tau) \in S_N$ the terms $\gamma^3 p (|\sigma|^2+|\tau|^2)$ and $\gamma^4 p |\sigma \tau|^2$ are of order $o(\frac 1 {Np})$.
Finally all $\eee^{\beta^2}$- and $\eee^{\frac{\beta^4}{N^2p^3}}$-terms are at most of constant order and appear identically in
both, $\E[T(\sigma)T(\tau)]$ and $\E [T(\sigma)] \, \E [T(\tau)]$.
Thus we obtain in this case
\begin{align}
\sum_{(\sigma,\tau) \in S_N} & \Cov (T(\sigma), T(\tau))\notag \\
&=
e^{C_6+o(1)}\sum_{(\sigma,\tau) \in S_N} \eee^{\frac{\beta}2 (\frac{|\sigma|^2}N+\frac{|\tau|^2}N)}
\left(\eee^{\frac{\beta^2(1-p)}{4Np}(\frac{|\sigma\tau|^2}N-1)+o(\frac 1 {Np})}-\eee^{o\left(\frac{1}{Np}\right)}\right)
\notag \\
&=  e^{C_6+o(1)}\sum_{(\sigma,\tau) \in S_N} \eee^{\frac{\beta}2 (\frac{|\sigma|^2}N+\frac{|\tau|^2}N)}
\left(\frac{\beta^2(1-p)}{4Np}\left(\frac{|\sigma\tau|^2}N-1\right)+o\left(\frac {1} {Np}\right)
\right)\notag\\
&= e^{C_6+o(1)}\sum_{(\sigma,\tau) \in S_N} \eee^{\frac{\beta}2 (\frac{|\sigma|^2}N+\frac{|\tau|^2}N)}
\left(\frac{\beta^2(1-p)}{4Np}\left(\frac{|\sigma\tau|^2}N-1\right)\right)+o\left(\frac {4^N} {Np}\right)\label{est_cov_2}
\end{align}
for a constant $C_6=C_6(\beta)$.
To obtain the third equality we expanded
$$
\eee^{o(\frac{1}{Np})}=1+o\left(\frac{1}{Np}\right)
$$
and
$$
\eee^{\frac{\beta^2(1-p)}{4Np}(\frac{|\sigma\tau|^2}N-1)+o(\frac 1 {Np})}= 1 + \frac{\beta^2(1-p)}{4Np}\left(\frac{|\sigma\tau|^2}N-1\right)+o\left(\frac 1 {Np}\right).
$$
Indeed, note that the second order term is
$$
\frac{1}{2} \left( \frac{\beta^2(1-p)}{4Np}\left(\frac{|\sigma\tau|^2}N-1\right)+ o\left(\frac 1 {Np}\right)\right)^2 \le
C \left(\frac{N (Np)^{\frac 15}}{N Np}\right)^2= \mathcal{O}\left((Np)^{-\frac 85}\right) =o\left(\frac 1{Np}\right)
$$
by construction of $S_N$.

Next we want to extend the range of summation in \eqref{est_cov_2} to all of $\{\pm 1\}^N \times \{\pm 1\}^N $. To this end, as in the first part of the proof set
$$
D_N\coloneqq\{(x,y,z)\in\R^3\colon |x|>(Np)^{\frac 1{10}} \text{ or } |y|>(Np)^{\frac 1{10}} \text{ or } |z|>(Np)^{\frac 1{10}} \}
$$
and use \eqref{LCLT} to estimate
\begin{multline*}
\left|\sum_{(\sigma,\tau) \notin S_N} \eee^{\frac{\beta}2 (\frac{|\sigma|^2}N+\frac{|\tau|^2}N)}
\left(\frac{\beta^2(1-p)}{4Np}\left(\frac{|\sigma\tau|^2}N-1\right)\right)\right|
\\
\le
C_1 4^{N} \frac{\beta^2}{4Np}N^{-3/2}
\sum_{\substack{(k,l,m)\in\Z^3\\N^{-1/2}(k,l,m)\in D_N}} \eee^{-\delta \frac{k^2+l^2 + m^2} {2N}}
\left|\frac{m^2}N-1\right|
\end{multline*}
for some $\delta >0$ small enough.
Similar to what we did in the first part of the proof
\begin{equation*}
N^{-3/2}
\sum_{\substack{(k,l,m)\in\Z^3\\N^{-1/2}(k,l,m)\in D_N}} \eee^{-\delta \frac{k^2+l^2 + m^2} {2N}}
\left|\frac{m^2}N-1\right| \le
\int_{\frac 13 D_N} \eee^{-\delta \frac{x^2+y^2 + z^2} {2}} |z^2-1|\dd x\, \dd y \,  \dd z
\end{equation*}
and the right hand side goes to zero as $N \to \infty$. Thus
$$
\sum_{(\sigma,\tau) \notin S_N} \eee^{\frac{\beta}2 (\frac{|\sigma|^2}N+\frac{|\tau|^2}N)}
\left(\frac{\beta^2(1-p)}{4Np}\left(\frac{|\sigma\tau|^2}N-1\right)\right)= o\left(\frac{4^N}{Np}\right).
$$
Thus, we obtain 
\begin{multline*}
\sum_{(\sigma,\tau) \in S_N}  \Cov (T(\sigma), T(\tau)) \\
=
e^{C_6+o(1)}\sum_{(\sigma,\tau) \in \{\pm 1\}^N\times \{\pm 1\}^N} \eee^{\frac{\beta}2 (\frac{|\sigma|^2}N+\frac{|\tau|^2}N)}
\left(\frac{\beta^2(1-p)}{4Np}\left(\frac{|\sigma\tau|^2}N-1\right)
\right)+o\left(\frac{4^N}{Np}\right).
\end{multline*}

%
%
%
As observed above, if $\sigma$ and $\tau$ are taken independently and uniformly at random from $\{-1,+1\}^N$, the vectors $(\sigma_i,\tau_i,\sigma_i\tau_i)$, $1\leq i\leq N$, are i.i.d.\ and we may apply the three-dimensional Central Limit Theorem to show that
$(\frac{|\sigma|}{\sqrt N},\frac{|\tau|}{\sqrt N},\frac{|\sigma \tau|}{\sqrt N})$ converges to a Gaussian random vector with expectation vector $0$ and identity covariance matrix.
We thus have
\begin{multline*}
4^{-N}\sum_{(\sigma,\tau) \in \{\pm 1\}^N\times \{\pm 1\}^N} \eee^{\frac{\beta}2 (\frac{|\sigma|^2}N+\frac{|\tau|^2}N)}
\left(\frac{|\sigma\tau|^2}N-1\right)
\\
\to C_7 \int_{\mathbb R^3} \eee^{(\frac{\beta-1}2)(x^2+y^2)-\frac{z^2}2}(z^2-1) \, dx\, dy\, dz.
\end{multline*}
Indeed, this is true, using the same uniform integrability argument as in \eqref{eq:UI}.
But
$$
\int_{\mathbb R} \eee^{-\frac{z^2}2}(z^2-1)dz=0
$$
and thus
$$
\sum_{(\sigma,\tau) \in \{\pm 1\}^N\times \{\pm 1\}^N} \eee^{\frac{\beta}2 (\frac{|\sigma|^2}N+\frac{|\tau|^2}N)}
\left(\frac{\beta^2(1-p)}{4Np}\left(\frac{|\sigma\tau|^2}N-1\right)\right)=o\left(\frac{4^N}{Np}\right),
$$
which yields
$$
Np \sum_{(\sigma,\tau) \in S_N}
 \Cov (T(\sigma), T(\tau)) =o(4^N).
$$
Putting the estimates for typical and atypical spins together shows the assertion.
\end{proof}

\subsection{Completing the proof of Theorem~\ref{CLT_ZN_2}}
To finish the proof of our Theorem~\ref{CLT_ZN_2} let us introduce some notation. For a sequence of random variables $(X_N)_{N\in \mathbb N}$ and sequence of positive real numbers $(a_N)_{N \in \mathbb N}$ we will say that $X_N=o_{L^2}(a_N)$ if $X_N/a_N\to 0$ in $L^2$, that is if $\E(X^2_N)= o(a_N^2)$. In the sense of this notation, what we showed in Proposition \ref{var_sum} can be written as
$$
\frac{\hat Z_N(\beta)}{\E \hat Z_N(\beta)} - 1 = o_{L^2}\left(\frac 1 {\sqrt{Np}}\right).
$$
In view of~\eqref{eq:Z_N_hat_Z_N} this implies
\begin{align}
\frac{Z_N(\beta)}{\E \hat Z_N(\beta)}
&=
\frac{\hat Z_N(\beta)}{\E \hat Z_N(\beta)}
\exp\left( \alpha_N \xi_N + \beta_N\eta_N\right)\notag\\
&=
\left(1 + o_{L^2}\left(\frac 1 {\sqrt{Np}}\right)\right)\exp\left( \alpha_N \xi_N + \beta_N\eta_N\right),
\label{prep_proof_th1.2}
\end{align}
where we defined two sequences of constants
\begin{equation}\label{eq:alpha_N_beta_N}
\alpha_N := \frac{\beta\sqrt{1-p}}{2\sqrt {Np}}
\;\;\;
\text{ and }
\;\;\;
\beta_N:= \frac{\beta^2 (1-2p)\sqrt{1-p}}{8N p^{3/2}}.
\end{equation}
Introduce the centered variables
\begin{equation}\label{eq:centered}
\overline{\xi_N} \coloneqq \frac 1 {\sqrt{N p(1-p)}} \sum_{i=1}^N \overline{\vep_{i,i}},
\;\;\;
\text{ and }
\;\;\;
\overline{\eta_N} \coloneqq \frac 1 {N \sqrt{p(1-p)}} \sum_{i,j=1}^N \overline{\vep_{i,j}},
\end{equation}
where
$
\overline{\vep_{i,j}}\coloneqq \vep_{i,j}-p
$
is the centered version of $\vep_{i,j}$.
With this notation, we can rewrite~\eqref{prep_proof_th1.2} as follows:
\begin{equation}\label{prep_proof_th1.2_rewr}
 \frac{Z_N(\beta)}{e^{\frac \beta 4 + \frac{\beta}{8p}} \E \hat Z_N(\beta)}
=
\left(1 + o_{L^2}\left(\frac 1 {\sqrt{Np}}\right)\right)
\exp\left( \alpha_N \overline{\xi_N} + \beta_N\overline{\eta_N}\right),
\end{equation}
The idea is now to replace the exponential on the right-hand side by its linearization which, as the next lemma states, has asymptotically Gaussian fluctuations.
\begin{lemma}\label{lem:clt_for_theta_N}
Assuming that $p^2 N  \to c^2 \in [0,\infty)$ (where $c$ may be $0$) and $p^3 N^2\to\infty$, we have
$$
Np^{3/2} \left(\alpha_N \overline{\xi_N} + \beta_N\overline{\eta_N}\right) \to \mathcal{N}\left(0,\frac{\beta^4}{64} + \frac{\beta^2c}{4}\right)
$$
in distribution.
\end{lemma}
\begin{proof}
Let us write $\theta_N := \alpha_N \overline{\xi_N} + \beta_N\overline{\eta_N}$ to shorten the notation. Using~\eqref{eq:centered} and~\eqref{eq:alpha_N_beta_N}, we have
\begin{equation}\label{eq:tech1}
Np^{3/2} \theta_N
=
\sum_{i\neq j} \frac{\beta^2(1-2p)}{8N\sqrt p} \overline{\vep_{i,j}}
+
\sum_{i=1}^N
\left(\frac{\beta^2(1-2p)}{8N\sqrt p} + \frac {\beta\sqrt p}{2}\right) \overline{\vep_{i,i}}.
\end{equation}
The right-hand side is a sum of independent, zero-mean random variables. The variance of the right-hand side is
\begin{align}
\V (Np^{3/2} \theta_N)
&=
\sum_{i\neq j} \frac{\beta^4(1-2p)^2 (1-p)}{64N^2}  + \sum_{i=1}^N
\left(\frac{\beta^2(1-2p)}{8N\sqrt p} + \frac {\beta\sqrt p}{2}\right)^2 p(1-p)\notag
\\
&\to
\frac{\beta^4}{64} + \frac{\beta^2c}{4}\label{eq:var_asympt}
\end{align}
by our assumptions. Moreover, since $1/(N\sqrt p)= o(\sqrt p)  \to 0$, all random variables on the right-hand side are uniformly bounded by some $\delta_N\to 0$. It follows that the Lindeberg CLT can be applied, which yields the assertion.
\end{proof}

\begin{lemma}\label{lem:linearization}
We have
$$
\exp\left( \alpha_N \overline{\xi_N} + \beta_N\overline{\eta_N}\right)
=
1 + \alpha_N \overline{\xi_N} + \beta_N\overline{\eta_N} + o_{L^2}\left(\frac{1}{Np^{3/2}}\right).
$$
\end{lemma}
\begin{proof}
Recall the notation $\theta_N= \alpha_N \overline{\xi_N} + \beta_N\overline{\eta_N}$.  Using the Taylor series for the exponential function and writing the remainder term in the Lagrange form, we have
$$
|e^{\theta_N} - 1 - \theta_N| = \frac 12 e^{r_N} r_N^2 \leq \frac 12 e^{|\theta|_N} \theta_N^2,
$$
for some $r_N$ with $|r_N|\leq |\theta_N|$. It follows that
$$
N^2 p^3 \E \left(|e^{\theta_N} - 1 - \theta_N|^2\right) \leq \frac 14 N^2 p^3 \E (e^{2|\theta|_N} \theta_N^4).
$$
Our aim is to prove that the right-hand side goes to $0$. By the Cauchy-Schwarz inequality, it suffices to show that
$$
\E (e^{4|\theta|_N}) =\mathcal O(1)
\quad \text{ and } \quad
N^4 p^6 \E (\theta_N^8) =o(1).
$$
First of all note that $N^4 p^6 \theta_N^8 = (Np^{3/2}\theta_N)^8 / (Np^{3/2})^4 \to 0$ in distribution because the numerator converges in distribution to the fourth power of a normal variable by Lemma~\ref{lem:clt_for_theta_N}. What is missing is the corresponding convergence of expectations, for which we need the uniform integrability.  Thus, it suffices to check that
\begin{equation}\label{eq:uniform_integr}
\E (e^{4|\theta|_N}) =\mathcal O(1)
\quad \text{ and } \quad
\E ((Np^{3/2}\theta_N)^{16}) = \mathcal O(1).
\end{equation}
To prove both assertions, we use the Bernstein inequality; see Theorem 2.8.4 in~\cite{vershynin_book}. It states that for independent, zero-mean random variables $X_1,\ldots,X_m$ with $|X_i|<K$ we have
$$
\P\left(\left|\sum_{i=1}^m X_i\right| \geq t \right) \leq 2 \exp\left(-\frac {t^2/2}{\V \sum_{i=1}^m X_i + Kt/3}\right)
$$
for all $t\geq 0$.
Applying this to the independent random variables on the right-hand side of~\eqref{eq:tech1} with $K= 3C_1\sqrt p$, we obtain
$$
\P\left(\left|Np^{3/2} \theta_N\right| \geq t \right) \leq 2 \exp\left(-\frac {t^2/2}{C_2 + C_1 t \sqrt p}\right)
\leq C_3 \exp\left(-C_4 t\right)
$$
for all $t>0$. Here, $C_1,\ldots, C_4$ are some positive constants and we used~\eqref{eq:var_asympt} to prove the estimate $\V \sum_{i=1}^m X_i\leq C_2$.   This exponential tail estimate is uniform in $N$ and immediately implies the second claim in~\eqref{eq:uniform_integr}. To prove the first claim, recall that $Np^{3/2}\to\infty$, hence for sufficiently large $N$ we have
$
\P\left(\left|\theta_N\right| \geq t \right) \leq C_3 \exp (- t),
$
which yields the first claim in~\eqref{eq:uniform_integr}.
\end{proof}

We are now in position to complete the proof of Theorem~\ref{CLT_ZN_2}. Assume that $p^2N\to c^2\in [0,\infty)$ and $p^3 N^2\to\infty$, where $c=0$ is possible. Using~\eqref{prep_proof_th1.2_rewr} and Lemma~\ref{lem:linearization}, we can write
\begin{align}
\frac{Z_N(\beta)}{e^{\frac \beta 4 + \frac{\beta}{8p}} \E \hat Z_N(\beta)}
&=
\left(1 + o_{L^2}\left(\frac 1 {\sqrt{Np}}\right)\right)
\left(1 + \alpha_N \overline{\xi_N} + \beta_N\overline{\eta_N} + o_{L^2}\left(\frac{1}{Np^{3/2}}\right)\right)\notag\\
&=
1 + \alpha_N \overline{\xi_N} + \beta_N\overline{\eta_N} + o_{L^2}\left(\frac 1 {Np^{3/2}}\right).\label{eq:L2_tech}
\end{align}
Here, we used that $\alpha_N \overline{\xi_N} + \beta_N\overline{\eta_N} = o_{L^2}(1)$ (which follows from~\eqref{eq:var_asympt}) as well as the rule $o_{L^2}(a_N) o_{L^2}(b_N) = o_{L^2}(a_nb_N)$ which easily follows from the Cauchy-Schwarz inequality. Multiplying by $Np^{3/2}$, we obtain
$$
N p^{3/2} \left(\frac{Z_N(\beta)}{e^{\frac \beta 4 + \frac{\beta}{8p}} \E \hat Z_N(\beta)} - 1\right)
=
Np^{3/2}(\alpha_N \overline{\xi_N} + \beta_N\overline{\eta_N}) + o_{L^2}(1).
$$
Recalling the assumptions  $p^2N\to c^2\in [0,\infty)$ and $p^3 N^2\to\infty$, we can apply Lemma~\ref{lem:clt_for_theta_N} to obtain
\begin{equation}\label{eq:CLT_wrong_norming}
N p^{3/2} \left( \frac{Z_N(\beta)}{e^{\frac \beta 4 + \frac{\beta}{8p}} \E \hat Z_N(\beta)} - 1\right)
\to
\mathcal{N}\left(0,\frac{\beta^4}{64} + \frac{\beta^2c}{4}\right)
\end{equation}
in distribution. Up to normalization sequences, this coincides with the claim of Theorem \ref{CLT_ZN_2}. 

To complete the proof, we need essentially to replace $e^{\frac \beta 4 + \frac{\beta}{8p}} \E \hat Z_N(\beta)$ by $\E Z_N(\beta)$ in the above statement. This can be justified as follows. Multiplying~\eqref{eq:L2_tech} by $\sqrt{Np}$, using the fact that $L^2$-convergence implies convergence of expectations, and recalling that the expectation of $\alpha_N \overline{\xi_N} + \beta_N\overline{\eta_N}$ vanishes, we obtain
\begin{equation}\label{eq:eps_N}
\eps_N := Np^{3/2} \left( \frac{\E Z_N(\beta)}{e^{\frac \beta 4 + \frac{\beta}{8p}} \E \hat Z_N(\beta)} -1\right) = o(1).
\end{equation}
Using this notation, we can write
\begin{align*}
N p^{3/2} \left(\frac{Z_N(\beta)}{\E Z_N(\beta)} - 1\right)
&=
N p^{3/2} \left(Z_N(\beta)\cdot \frac{1+\eps_N N^{-1}p^{-3/2}}{\E Z_N(\beta)} - 1\right) - \eps_N \frac{Z_N(\beta)}{\E Z_N(\beta)}\\
&=
N p^{3/2} \left(\frac{Z_N(\beta)}{e^{\frac \beta 4 + \frac{\beta}{8p}} \E \hat Z_N(\beta)} - 1\right) - \eps_N \frac{Z_N(\beta)}{\E Z_N(\beta)}\\
&\to
\mathcal{N}\left(0,\frac{\beta^4}{64} + \frac{\beta^2c}{4}\right)
\end{align*}
in distribution, where in the last line we used~\eqref{eq:CLT_wrong_norming} and the fact that $\eps_N \frac{Z_N(\beta)}{\E Z_N(\beta)} \to 0$ in distribution and even in $L^1$ because $\eps_N\to 0$ by~\eqref{eq:eps_N}.
This proves the second claim of Theorem~\ref{CLT_ZN_2}. To prove the first claim, note that its assumption $p\sqrt N \to c\in (0,\infty)$ implies that $Np^{\frac{3}{2}}\sim c^2 p^{-\frac 12}\sim c^{3/2}N^{1/4}$.
\hfill $\Box$

\begin{rem}
 The above proof could be modified to show the validity of Theorem~\ref{CLT_ZN_1}. However, we prefer to give the proof that we outlined in Section 3 because the technique used there may be of interest on its own right.
\end{rem}

\section{Proof of Theorem \ref{CLT_ZN_3} and Theorem~\ref{CLT_ZN_4}}
For the proof of both theorems we will need the following modification of the partition function that was already used in \cite{KLS19b}.
For fixed $0<\beta <1$ we define $\gamma := \frac{\beta}{2 Np}$ and
\begin{equation}\label{eq:tildeT}
\tilde T(\sigma) \coloneqq \exp\left(\gamma \sum_{i,j=1}^N \vep_{i,j} \sigma_i \sigma_j - \log \cosh(\gamma)\sum_{i,j=1}^N \vep_{i,j}\right).
\end{equation}
Moreover, let
\begin{equation}\label{tildeZN}
\tilde Z_N(\beta)\coloneqq \sum_{\sigma \in \{-1, +1\}^N} \tilde T(\sigma)
=\exp\left(-\log \cosh(\gamma) \sum_{i,j=1}^N \vep_{i,j}\right)  Z_N(\beta).
\end{equation}
In \cite{KLS19b}, Proposition 3.5 and formula (3.15), we have shown the following lemma playing the same role as Proposition~\ref{var_sum} in the previous section.
\begin{lemma}\label{lem:tildeZN}
In the regime where $Np \to \infty$ and $0< \beta <1$ we have that
\begin{equation}\label{eq:in_probab}
\frac{ \tilde Z_N(\beta)}{\E \tilde Z_N(\beta)} \to 1
\end{equation}
in $L^2$ and hence in probability.
\end{lemma}
\begin{proof}[Proof of Theorems \ref{CLT_ZN_3} and~\ref{CLT_ZN_4}]
Let us center the variables $\vep_{i,j}$ and introduce
$$
\overline{\vep_{i,j}}\coloneqq \vep_{i,j}-p.
$$
Then,
$$
\tilde Z_N(\beta)= \exp\left(-\log \cosh(\gamma) \sum_{i,j=1}^N \overline{\vep_{i,j}}\right) \eee^{-\log\cosh(\gamma)N^2p} Z_N(\beta).
$$
After taking the logarithm, Lemma \ref{lem:tildeZN} therefore implies that
\begin{equation}\label{eq:proof_log_fluct_1}
\left(\log Z_N(\beta)-\log \E \tilde Z_N(\beta)-\log\cosh(\gamma)N^2p\right) - \log \cosh(\gamma) \sum_{i,j=1}^N \overline{\vep_{i,j}}\to 0
\end{equation}
in probability. On the other hand,
\begin{equation}\label{eq:proof_log_fluct_2}
\frac 1{\sqrt{N^2p (1-p)}} \sum_{i,j} \overline{\vep_{i,j}} \to
\mathcal{N}(0,1)
\end{equation}
in distribution, by the central limit theorem.

We now assume that  $N^2 p^3 \to c \in (0, \infty)$, as in Theorem~\ref{CLT_ZN_3}. This also implies that $p \to 0$ and  we have
$$
\log \cosh{(\gamma)}\sqrt{N^2p (1-p)}= \left(\frac{\gamma^2}{2}+\mathcal{O}(\gamma^4)\right)\sqrt{N^2p (1-p)}=
\frac{\beta^2}{8 \sqrt c}+o(1).
$$
Therefore, \eqref{eq:proof_log_fluct_1} and~\eqref{eq:proof_log_fluct_2} yield
$$
\log Z_N(\beta)-\log \E \tilde Z_N(\beta)-\log\cosh(\gamma)N^2p \to \mathcal{N}\left(0,\frac{\beta^4}{64c}\right)
$$
in distribution. According to Proposition~3.4 in~\cite{KLS19b} with $g=1$,
\begin{equation}\label{eq:prop_3.4}
\log \E \tilde Z_N(\beta) = N\log 2 - \frac {\beta^2}{8} - \frac 12 \log (1-\beta) + o(1).
\end{equation}
Moreover,
$$
\log\cosh(\gamma)N^2 p  = N^2p \left(\frac{\beta^2}{8N^2p^2} - \frac {\beta^4+o(1)}{192 N^4p^4}\right) = \frac{\beta^2}{8p} - \frac{\beta^4}{192 c} + o(1).
$$
Taking everything together we arrive at the assertion of Theorem~\ref{CLT_ZN_3}.

Let us now assume that $p^3 N^2 \to 0$ as in Theorem~\ref{CLT_ZN_4}. Then,
\begin{equation}\label{eq:proof_log_fluct_3}
\log \cosh{(\gamma)}\sqrt{N^2p (1-p)}= \left(\frac{\gamma^2+o(1)}{2}\right)\sqrt{N^2p (1-p)}\sim
\frac{\beta^2}{8 Np^{3/2}}\to\infty.
\end{equation}
Inserting~\eqref{eq:prop_3.4} (which holds under $Np\to\infty$) into~\eqref{eq:proof_log_fluct_1} and dividing everything by $\log \cosh{(\gamma)}\sqrt{N^2p (1-p)}\to\infty$, we arrive at
$$
\frac{\log Z_N(\beta)
-
N\log 2
-\log\cosh(\gamma)N^2p}{\log \cosh{(\gamma)}\sqrt{N^2p (1-p)}}- \frac {\sum_{i,j=1}^N \overline{\vep_{i,j}}} {\sqrt{N^2p (1-p)}} \to 0
$$
in probability. Recalling~\eqref{eq:proof_log_fluct_2} and~\eqref{eq:proof_log_fluct_3}, we arrive at the assertion of Theorem~\ref{CLT_ZN_4}.
\end{proof}


\end{document}